\documentclass[11pt]{article}

\usepackage[T1,T2A]{fontenc}
\usepackage[utf8]{inputenc}
\usepackage{amsfonts}
\usepackage{amssymb}
\usepackage{epsfig}
\usepackage{tikz}

\usepackage{theorem}

\usepackage{color}

\newtheorem{theorem}{Theorem}[section]
\newtheorem{lemma}[theorem]{Lemma}
\newtheorem{corollary}[theorem]{Corollary}

{\theorembodyfont{\rmfamily} % No need in ``\rm'' command after ``\begin''

\newtheorem{remark}[theorem]{Remark}

}

\def\Q{\mathbb{Q}}

\def\A{\mathcal{A}}
\def\B{\mathcal{B}}
\def\C{\mathcal{C}}
\def\D{\mathcal{D}}
\def\D{\mathcal{D}}
\def\E{\mathcal{E}}
\def\F{\mathcal{F}}
\def\G{\mathcal{G}}
\def\H{\mathcal{H}}
\def\I{\mathcal{I}}
\def\J{\mathcal{J}} 
\def\K{\mathcal{K}}
\def\L{\mathcal{L}} 
\def\M{\mathcal{M}}
\def\N{\mathcal{N}}
\def\R{\mathcal{R}}
\def\T{\mathcal{T}}
\def\U{\mathcal{U}}
\def\V{\mathcal{V}}
\def\W{\mathcal{W}}
\def\X{\mathcal{X}}
\def\Y{\mathcal{Y}}
\def\Z{\mathcal{Z}}

\title{Doubly Hurwitz Beauville groups}
\date{\today}
\author{Gareth A. Jones and Emilio Pierro}

\begin{document}

\maketitle

\begin{abstract}
If $\mathcal S$ is a Beauville surface $({\mathcal C}_1\times{\mathcal C}_2)/G$, then the Hurwitz bound implies that $|G|\le 1764\,\chi({\mathcal S})$, with equality if and only if the Beauville group $G$ acts as a Hurwitz group on both curves $\C_i$. Equivalently, $G$ has two generating triples of type $(2,3,7)$, such that no generator in one triple is conjugate to a power of a generator in the other. We show that this property is satisfied by alternating groups $A_n$, their double covers $2.A_n$, and special linear groups $SL_n(q)$ if $n$ is sufficiently large, but by no sporadic simple groups or simple groups $L_n(q)$ ($n\le 7$), ${}^2G_2(3^e)$, ${}^2F_4(2^e)$, ${}^2F_4(2)'$, $G_2(q)$ or ${}^3D_4(q)$ of small Lie rank.
\end{abstract}


\section{Introduction}

A {\em Beauville surface\/} is a complex algebraic surface of general type, of the form
\[{\mathcal S}=({\mathcal C}_1\times{\mathcal C}_2)/G\]
where ${\mathcal C}_1$ and ${\mathcal C}_2$ are the algebraic curves underlying regular dessins $\R_1$ and $\R_2$ of genera $g_1, g_2>1$, which have the same automorphism $G$ which acts freely on their product. (This construction was introduced by Beauville in~\cite{Bea}; for further details see the surveys~\cite{Fairbairn-15, Jones-14} or~\cite[Chapter 11]{JW-16}.) A finite group $G$ arises in this way (and is then called a {\em Beauville group}) if and only if it has generating triples $(x_i, y_i, z_i)$ of hyperbolic types for $i=1, 2$, such that no non-identity element of $G$ has fixed points on both curves, or equivalently no non-identity power of $x_1, y_1$ or $z_1$ is conjugate in $G$ to a power of $x_2, y_2$ or $z_2$.

In this situation, $\mathcal S$ has Euler characteristic
\[\chi({\mathcal S})=\frac{\chi({\mathcal C}_1)\chi({\mathcal C}_2)}{|G|}
=\frac{4(g_1-1)(g_2-1)}{|G|}.\]
Using the Hurwitz bound $|G|\le 84(g_i-1)$ (see~\cite{Hurwitz-1893}) we therefore have
\begin{equation}\label{1764ineq}
|G|\le 1764\,\chi({\mathcal S}),
\end{equation}
with equality if and only if $G$ acts on each curve $\C_i$ as a Hurwitz group, that is, $\C_i$ is uniformised by a normal subgroup of the triangle group
\[\Delta=\Delta(2,3,7)=\langle X,Y,Z\mid X^2=Y^3=Z^7=XYZ=1\rangle,\]
with quotient group $G$.  Alexander Zvonkin~\cite{Zvonkin-16} has raised the question of whether there exist such groups $G$ which act as Hurwitz groups in two essentially different ways, that is, which have two generating triples $(x_i, y_i, z_i)$ of type $(2,3,7)$ such that  no non-identity power of $x_1, y_1$ or $z_1$ is conjugate in $G$ to a power of $x_2, y_2$ or $z_2$, so that the bound in~(\ref{1764ineq}) is attained. We will call such a group a {\em doubly Hurwitz Beauville group}, or a {\em dHB group\/} for brevity. (By contrast we note that the first examples of Beauville surfaces~\cite{Bea}, where each $\mathcal C_i$ is the Fermat curve of degree $n$ coprime to $6$ and $G\cong C_n\times C_n$, so that $\chi(\mathcal S)=(n-3)^2$, come nowhere near attaining this bound.)

Each Hurwitz group is perfect, and is thus a cover of a non-abelian finite simple group, which is also a Hurwitz group. It is therefore of interest to consider various families of non-abelian finite simple groups, and more generally quasisimple groups, to determine which of their members are dHB groups. Here, being a Beauville group is no great restriction: it is easily seen that $A_5$ is not a Beauville group, and in 2005 Bauer, Catanese and Grunewald~\cite{BCG} conjectured that every other non-abelian finite simple group is a Beauville group. This was proved, first with finitely many possible exceptions, by Garion, Larsen and Lubotsky~\cite{GLL}, and then completely by Guralnick and Malle~\cite{GM}, a result which Fairbairn, Magaard and Parker~\cite{FMP} extended to all quasisimple groups except $A_5$ and its double cover $\tilde A_5\cong SL_2(5)$. However, being a dHB group is a much stronger condition than being both a Hurwitz group and a Beauville group.

As a perhaps surprising example, none of the 26 sporadic simple groups is a dHB group (see Theorem~\ref{sporadic}), even though twelve of them are Hurwitz groups: not even the Monster, shown by Wilson~\cite{Wilson-01} to be a Hurwitz group in many different ways, has a pair of generating triples of type $(2,3,7)$ which satisfy the Beauville condition.

As shown by Macbeath~\cite{Macbeath-69} the groups $G=L_2(p)=PSL_2(p)$, for primes $p\equiv\pm 1$ mod 7, are Hurwitz groups in three ways, since ${\rm Aut}\,G$ has three orbits on generating triples of type $(2,3,7)$. However, these groups have only one conjugacy class of involutions, so the Beauville non-conjugacy condition cannot be satisfied (a similar problem arises with their subgroups of orders 3 and 7). On the other hand, sufficiently large alternating groups $A_n$ do not suffer from these restrictions, and indeed they have been shown by Fuertes and Gonz\'alez-Diez~\cite{Fuertes-GD-10} to be Beauville groups for all $n\ge 6$, and by Conder~\cite{Conder-80} to be Hurwitz groups for all $n\ge 168$ (as well as for certain values $n<168$). We will adapt Conder's proof of this to establish the following theorem:

\begin{theorem}\label{maintheorem}
For each $r=0, 1, \ldots 13$ there exists an integer $N_r$, given in Table~\ref{Nr}, such that if $n\equiv r$ {\rm mod}~$(14)$ and $n\ge N_r$ then the alternating group $A_n$ is a doubly Hurwitz Beauville group. In particular, $A_n$ is a doubly Hurwitz Beauville group for all $n\ge 589$.
\end{theorem}

\begin{table}[ht]
\centering
\begin{tabular}{| p{0.35cm} | p{0.43cm} | p{0.43cm} | p{0.43cm} | p{0.43cm} | p{0.43cm} |
p{0.43cm} | p{0.43cm} | p{0.43cm} | p{0.43cm} | p{0.43cm} | p{0.43cm} |
p{0.43cm} | p{0.43cm} | p{0.43cm} |}
\hline
$r$ & 0 &1 & 2 & 3 & 4 & 5 & 6 & 7 & 8 & 9 & 10 & 11 & 12 & 13\\
\hline
$N_r$ & 294 & 589 & 394 & 367 & 396 & 439 & 510 & 329 & 540 & 457 & 430 & 459 & 432 & 447 \\
\hline
\end{tabular}
\caption{Values of $N_r$.}
\label{Nr}
\end{table}

The proof, in Section~\ref{altgps}, proceeds in two stages: first the construction in $A_n$ of pairs of triples of type $(2,3,7)$ which satisfy the Beauville condition, and then the proof that they both generate the whole group. The first stage is not so difficult, and can be performed within alternating groups of considerably lower degrees. Heuristically, if $n$ is sufficiently large one might expect almost all such triples to generate $A_n$, by Dixon's theorem~\cite{Dixon-69} on random generation of alternating and symmetric groups. However, {\em proving\/} that a given triple generates $A_n$ is more difficult, requiring the construction of a permutation with a suitable cycle of prime length, so that a well-known theorem of Jordan can be applied (see Section~\ref{techniques}). Having to do this simultaneously for two essentially different triples makes the task even more delicate, and explains the rather large lower bounds on $n$ in Table~\ref{Nr}. No doubt a more careful construction of suitable triples, perhaps not so heavily reliant on the hard work already done by Conder in~\cite{Conder-80}, could achieve more realistic lower bounds. Indeed, minor modifications to our proof yield a few smaller examples, such as $A_{246}$, while Conder~\cite{Conder-17} has recently used Magma to show that $A_{168}$ is a dHB group. (This is, in fact, the smallest example of a dHB group currently known to us.) However, one should not expect major improvements, since $A_n$ is never a dHB group if $n<168$ (see Theorem~\ref{168thm}).

In 2010, Pellegrini and Tamburini~\cite{PT} showed that the double cover $\tilde A_n$ of $A_n$ is a Hurwitz group for all $n\ge 231$; in Section~\ref{dblcov}, a small modification of the proof of Theorem~\ref{maintheorem} gives the following, where $N_r$ is as in Table~\ref{Nr}:

\begin{corollary}\label{covercor}
If $n\equiv r$ {\rm mod}~$(14)$ and $n\ge N_r+112$ then the double cover $\tilde A_n$ of $A_n$ is a doubly Hurwitz Beauville group. In particular, this applies for all $n\ge 701$.
\end{corollary}

In 2000, Lucchini, Tamburini and Wilson~\cite{LTW} showed that if $n\ge 287$ then $SL_n(q)$ is a Hurwitz group for each prime power $q$, and it immediately follows that the same applies to the simple group $L_n(q)=PSL_n(q)$. By combining our method of proof with theirs (an adaptation of Conder's in~\cite{Conder-80}, using permutation matrices in place of permutations), we will prove the following in Section~\ref{SLn}:

\begin{corollary}\label{SLnqcor}
If $n\equiv r$ {\rm mod}~$(14)$ and $n\ge N_r+42$ then the groups $SL_n(q)$ and $L_n(q)$ are doubly Hurwitz Beauville groups for all prime powers $q$. In particular, this applies for all $n\ge 631$.
\end{corollary}

In contrast with the above results, we will show in Section~\ref{simple} that there are no dHB groups among various other families of simple groups, either by arguments similar to those which apply to Macbeath's Hurwitz groups, or by a more detailed analysis of their generating triples of type $(2,3,7)$:

\begin{theorem}\label{simplethm}
If $G$ is a simple group of one of the following types,
\begin{enumerate}
\item a sporadic simple group,
\item $L_n(q)$ ($n\le 7$), ${}^2G_2(3^e)$, ${}^2F_4(2^e)$, ${}^2F_4(2)'$, $G_2(q)$ or ${}^3D_4(q)$,
\end{enumerate}
then $G$ is not a doubly Hurwitz Beauville group:
\end{theorem}

The groups listed here include many of the families of simple Hurwitz groups of small Lie rank.
Of course, several families remain, and there is a vast no-man's-land between these groups and those in Corollary~\ref{SLnqcor}, in which results about dHB groups will no doubt be harder to obtain.

Although the techniques we use are almost entirely group-theoretic, in order to make this paper accessible to geometers and others we have given rather more background detail than experts in group theory might require. For a similar reason, in Section~\ref{dessins} we will discuss some connections between the method of proof of Theorem~\ref{maintheorem}, the algebraic theory of maps on surfaces developed by Singerman and the first author in~\cite{JS-78}, and Grothendieck's theory of dessins d'enfants~\cite{Gro}.

\medskip

\noindent{\bf Acknowledgements} The second author was supported by the grant SFB 701 at Universit\"at Bielefeld and is grateful to Alastair Litterick for many helpful conversations. 

%%%%%%%%%%%%%%%

\section{Techniques}\label{techniques}

\subsection{Some classical results}\label{classical}

In the proof of Theorem~\ref{maintheorem}, we will follow Conder in applying the following theorem of Jordan (see~\cite[Theorem 13.9]{Wielandt-64}):

\begin{theorem}[Jordan]\label{Jordantheorem}
If $G$ is a primitive permutation group of degree $n$, containing a cycle of prime length with at least three fixed points, then $G\ge {\rm A}_n$.
\end{theorem}

The following corollary is a minor generalisation of a result of Everitt~\cite[Lemma~2.3]{Everitt-00}:

\begin{corollary}\label{Jordancor}
Let $G$ be a permutation group of finite degree $n$, with a transitive subgroup $H=\langle h_1,\ldots, h_r\rangle$. Suppose that $G$ contains an element $g$ with a cycle $c$ such that
\begin{enumerate}
\item $c$ has prime length $p\le n-3$ which is coprime to the lengths of all other cycles of $g$;
\item for each $i=1,\ldots, r$, there is an element of $c$ with its image under $h_i$ also in $c$.
\end{enumerate}
Then $G\ge {\rm A}_n$.
\end{corollary}

\noindent{\sl Proof. } By (1), replacing $c$ by a suitable power allows us to assume that $g=c$, so $g$ fixes all points outside $c$. It is now sufficient to show that $G$ is primitive, since then Jordan's Theorem completes the proof.

Suppose that $G$ is imprimitive. For each block of imprimitivity $B$ meeting $c$, $|B\cap c|$ must divide $p$. If $|B\cap c|=1$ for all such $B$, then since $g$ permutes these blocks transitively yet fixes all points outside $c$, it follows that $|B|=1$, contradicting imprimitivity. Hence $|B\cap c|=p$, that is, $c$ is contained in a single block $B$. But now condition~(2) implies that each generator $h_i$ of $H$ maps $B$ to itself, and hence so does $H$, contradicting the transitivity of $H$. Thus $G$ is primitive, as required. \hfill$\square$

\medskip

The following result of Frobenius~\cite{Frobenius-1896} allows one to count the triples of a given type, such as $(2,3,7)$, in a given finite group:

\begin{theorem}[Frobenius]\label{th:frobenius}
Let $\X$, $\Y$ and $\Z$ be conjugacy classes in a finite group $G$. Then the number $n(\X,\Y,\Z)=n_G(\X,\Y,\Z)$ of solutions in $G$ of the equation $xyz=1$, where $x\in\X$, $y\in\Y$ and $z\in\Z$, is given by the formula
\[\frac{|\X|\cdot |\Y|\cdot |\Z|}{|G|}
\sum_{\chi}\frac{\chi(x)\chi(y)\chi(z)}{\chi(1)}, \]
where the sum is over all irreducible complex characters $\chi$ of $G$.
\end{theorem}

Of course, if such triples exist one still has to determine whether any of them generate the whole group. For this, knowledge of the subgroup lattice of $G$, or at least a list of its maximal subgroups, is important. 

The numbers $n_G(\X,\Y,\Z)$ given by Theorem~\ref{th:frobenius} are sometimes referred to as `structure constants', since they are closely related to the structure constants for the centre of the group algebra ${\mathbb C}G$. This has a basis consisting of the class sums $\underline{\X}=\sum\{x\mid x\in\X\}$ for the conjugacy classes $\X$ of $G$. The coefficient of $\underline{\Z}$ in $\underline{\X}\cdot\underline{\Y}$ is the number of pairs $(x,y)\in\X\times\Y$ such that $xy=z$ for some fixed $z\in\Z$, and this is given by
\[\frac{|\X|\cdot |\Y|}{|G|}\sum_{\chi}\frac{\chi(x)\chi(y)\chi(z^{-1})}{\chi(1)}
=\frac{|\X|\cdot |\Y|}{|G|}\sum_{\chi}\frac{\chi(x)\chi(y)\overline{\chi(z)}}{\chi(1)}
\;=\frac{n(\X,\Y,\Z^{-1})}{|\Z|}\]
where $\Z^{-1}$ denotes the conjugacy class $\{z^{-1}\mid z\in\Z\}$.
\medskip

The following is a special case of a more general theorem of Singerman~\cite{Singerman-70} about Fuchsian groups. If $M$ is a subgroup of finite index $n$ in $\Delta$, then $\Delta$ acts as a transitive group of degree $n$ on the cosets of $M$. If $\alpha$, $\beta$ and $\gamma$ are the numbers of fixed points of $X$, $Y$ and $Z$ in this action, then $M$ has signature
\[(g; 2^{[\alpha]}, 3^{[\beta]}, 7^{[\gamma]})\]
with $g$ (the genus of ${\mathbb H}/M$, where $\mathbb H$ is the hyperbolic plane) given by the Riemann-Hurwitz formula
\begin{equation}\label{RHeqn}
g=1+\frac{1}{84}\left(n-21\alpha-28\beta-36\gamma\right).
\end{equation}
Conversely, if $\Delta$ has a transitive permutation representation of degree $n$, then the stabilisers of points form a conjugacy class of subgroups of index $n$ with this signature.

%%%%%%%%%

\subsection{Conder's diagrams}

In~\cite{Conder-80}, Marston Conder constructed quotients of $\Delta$ by applying Graham Higman's technique of `sewing coset diagrams together'. These diagrams illustrate pairs of permutations $x$ and $y$ of order $2$ and $3$, with product of order $7$, which generate a transitive group. Each diagram is a coset diagram for the stabiliser of a point, with respect to the generators $x$ and $y$: small triangles, taken anticlockwise, represent 3-cycles of $y$, and heavy dots represent its fixed points, while long edges represent 2-cycles of $x$. An equivalent but more economical representation is to regard $x$ and $y$ as the canonical generators for a map on an oriented surface, as in~\cite{JS-78}, permuting arcs (directed edges) by reversing them or rotating them around their incident vertices: thus each of Conder's small triangles is shrunk to a point, so that vertices of valency 1 and 3 represent fixed points and 3-cycles of $y$, while edges and free edges represent 2-cycles and fixed points of $x$. The faces of the map then correspond to the cycles of $z:=(xy)^{-1}$; each face has seven sides or one, with a free edge counting as one side. The group $\langle x, y\rangle$, a quotient of $\Delta$ since $x^2=y^3=z^7=xyz=1$, is the cartographic group of the map, or the monodromy group of the corresponding dessin d'enfant~\cite{JW-16} in Grothendieck's terminology. (The connections with maps and dessins will be discussed in more detail in Section~\ref{dessins}.) For the rest of this paper, the terms `diagram' and `map', and the symbols $x, y$ and $z$, will always have the above meaning.

\begin{figure}[h!]
\begin{center}
\begin{tikzpicture}[scale=0.4, inner sep=0.7mm]

\draw [rounded corners, thick] (4,0) to (4,4) to (-4,4) to (-4,-4) to (4,-4) to (4,0);
\draw [thick] (1,-4) to (0,-2.5) to (-1,-4);
\draw [thick] (0,1) to (1,2.5) to (-1,2.5) to (0,1);
\draw [thick] (0,-2.5) to (0,1);
\draw [thick] (4,1) to (5.5,0) to (4,-1);
\draw [thick] (-4,1) to (-5.5,0) to (-4,-1);
\node (A) at (8,0) [shape=circle, fill=black] {};
\node (B) at (-8,0) [shape=circle, fill=black] {};
\draw [thick] (A) to (5.5,0);
\draw [thick] (B) to (-5.5,0);
\node at (6,-3) {$A$};

%%%%%%

\draw [thick] (17,0) arc (0:360:2);
\node (a) at (17,0) [shape=circle, fill=black] {};
\node (a') at (13,0) [shape=circle, fill=black] {};
\node (b) at (15,-2) [shape=circle, fill=black] {};
\node (c) at (15,-0.5) [shape=circle, fill=black] {};
\node (d) at (19.5,0) [shape=circle, fill=black] {};
\node (d') at (10.5,0) [shape=circle, fill=black] {};

\draw [thick] (b) to (c);
\draw [thick] (16,0.5) to (c) to (14,0.5);
\draw [thick] (a) to (d);
\draw [thick] (a') to (d');
\node at (18,-3) {$\A$};

\end{tikzpicture}

\end{center}
\caption{Conder's diagram $A$ and the corresponding map $\A$.}
\label{DiagramA}
\end{figure}
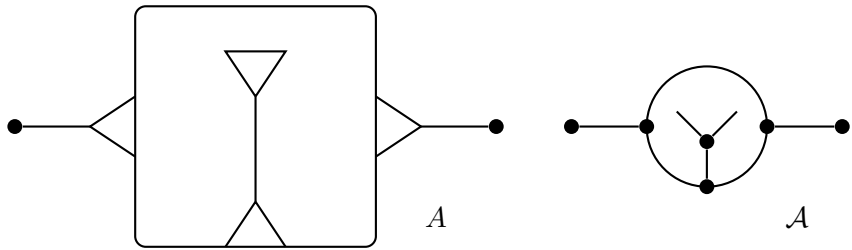

For example, Figure~\ref{DiagramA} shows a typical permutation diagram (called $A$ in~\cite{Conder-80}), with the corresponding map $\A$ (of genus 0). The permutations $x, y$ and $z$ are elements of the symmetric group $S_{14}$ with cycle structures $1^22^6$, $1^23^4$ and $7^2$. The cartographic group $\langle x, y\rangle$ is the simple group $L_2(13)$ in its natural representation of degree $14$ on the projective line ${\mathbb P}^1(\mathbb F_{13})$.

\medskip

\medskip

%%%%%%%%%%%%%%%%%%%

The diagrams Conder used, such as $A$ in Figure~\ref{DiagramA}, are all invariant under a reflection $t$ which inverts $x$ and $y$, so that $\langle x, y, t\rangle$ ie the {extended cartographic group of the map, a quotient of the extended triangle group
\[\Delta^*=\Delta[2,3,7]=\langle X,Y,T \mid X^2=Y^3=T^2=(XY)^7=(XT)^2=(YT)^2=1\rangle\]
which contains $\Delta$ with index $2$. Having odd order, $x$ and $z$ are always even permutations, and hence so is $y$, whereas $t$ could be even or odd.

\begin{table}[ht]
\centering
\begin{tabular}{| p{1.4cm} | p{1.1cm} | p{1cm} | p{2.1cm} | p{1.3cm} | p{3.1cm} |}
\hline
Diagram & Degree & Parity & Fixed points & Handles & Cycle-lengths of $w$ \\
\hline\hline
$A$ & $14$ & $+$ & $2,2, 0$ & $1, 0, 0$ & $1, 13$\\
\hline
$B$ & $15$ & $+$ & $3, 0, 1$ & $0, 2, 1$ & $3, {\bf 5}, 7$\\
\hline
$C$ & $21$ & $-$  & $5, 0, 0$ & $1, 0, 1$ & $1, {\bf 4}, {\bf 8}, 8$\\
\hline
$D$ & $22$ & $-$  & $2, 1, 1$ & $0, 1, 0$ & $5, 6, 11$\\
\hline
$E$ & $28$ & $+$ & $4, 1, 0$ & $1, 1, 0$ & $1, 9^3$\\
\hline
$F$ & $30$ & $+$ & $2, 0, 2$ & $0, 1, 0$ & $15^2$\\
\hline
$G$ & $42$ & $+$ & $6, 0, 0$ & $3, 0, 0$ & $1^3, 13^3$\\
\hline
$H$ & $42$ & $-$  & $6, 0, 0$ & $1, 0, 1$ & $1, 3, 10, 11, {\bf 17}$\\
\hline
$I$ & $57$ & $-$  & $5, 0, 1$ & $0, 2, 0$ & $4, 7, 8, 10, 13, 15$\\
\hline
$J$ & $72$ & $-$ & $4, 0, 2$ & $2, 0, 0$ & $1^2, 10, 11^2, 16, 22$\\
\hline
$K$ & $72$ & $+$ & $4, 0, 2$ & $1, 0, 0$ & $1, 5, {\bf 17}, 49$\\
\hline
$L$ & $102$ & $-$ & $2, 0, 4$ & $0, 1, 0$ & $21, 23, 58$\\
\hline
$M$ & $108$ & $+$ & $4, 0, 3$ & $1, 1, 0$ & $1, 12, 14, 19, 26, 36$\\
\hline
$N$ & $108$ & $+$ & $4, 0, 3$ & $1, 0, 1$ & $1, 9, 18, 20, 21, 39$\\
\hline

\end{tabular}
\caption{The 14 basic diagrams.}
\label{BasicDessins}
\end{table}

Conder constructed his diagrams from 14 basic diagrams $A,\ldots N$, drawn in~\cite{Conder-80}; the corresponding maps $\A,\dots, \N$ are shown in the Appendix. Table~\ref{BasicDessins} gives information about these diagrams which we will need later, namely the degree of the permutation group, the parity of $t$, the numbers of fixed points of $x, y$ and $z$, the numbers of $(k)$-handles for $k=1, 2$ and $3$ (see \S2.3), and the cycle structure of the permutation $w:=xyt$, with an entry $l^m$ indicating $m$ cycles of length $l$. (The significance of these cycle-lengths, together with the reason for showing some of them in bold face, will be

%%%%%%%%%%%%%%%%%%%%%%%%%%

\subsection{Composing diagrams and maps}

In~\cite{Conder-80}, Conder described a method for joining pairs of diagrams together to produce larger diagrams. Here we will reinterpret and illustrate this idea in terms of the corresponding maps. A $(k)$-{\em handle}, for $k=1, 2$ or $3$, is a pair $a, b$ of fixed points of the involution $x$ satisfying $b=a(xy)^k=at$ (see Figure~\ref{handles}). 

%%%%%%%%%%%%%%%%%%%%

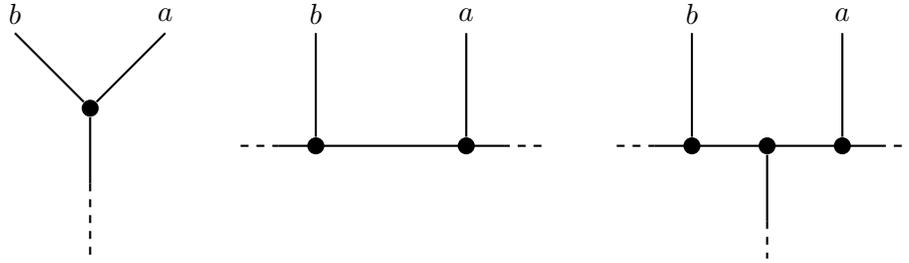
\begin{figure}[h!]
\begin{center}
\begin{tikzpicture}[scale=0.5, inner sep=0.8mm]

\node (1) at (-8,3) [shape=circle, fill=black] {};
\draw [thick] (1) to (-8,1);
\draw [thick, dashed] (-8,1) to (-8,-1);
\draw [thick] (-10,5) to (1) to (-6,5);
\node at (-6,5.5) {$a$};
\node at (-10,5.5) {$b$};

%%%%%%%

\node (2) at (-2,2) [shape=circle, fill=black] {};
\node (3) at (2,2) [shape=circle, fill=black] {};
\draw [thick] (2) to (-2,5);
\draw [thick] (3) to (2,5);
\draw [thick] (-3,2) to (3,2);
\draw [thick, dashed] (-4,2) to (-3,2);
\draw [thick, dashed] (4,2) to (3,2);
\node at (2,5.5) {$a$};
\node at (-2,5.5) {$b$};

%%%%%%

\node (4) at (8,2) [shape=circle, fill=black] {};
\node (5) at (12,2) [shape=circle, fill=black] {};
\node (6) at (10,2) [shape=circle, fill=black] {};
\draw [thick] (4) to (8,5);
\draw [thick] (5) to (12,5);
\draw [thick] (7,2) to (13,2);
\draw [thick, dashed] (6,2) to (7,2);
\draw [thick, dashed] (14,2) to (13,2);
\draw [thick] (6) to (10,0);
\draw [thick, dashed] (10,0) to (10,-1);
\node at (12,5.5) {$a$};
\node at (8,5.5) {$b$};

\end{tikzpicture}

\end{center}
\caption{$(k)$-handles for $k=1, 2, 3$.}\label{handles}
\end{figure}

%%%%%%%%%%%%%%%%%%%

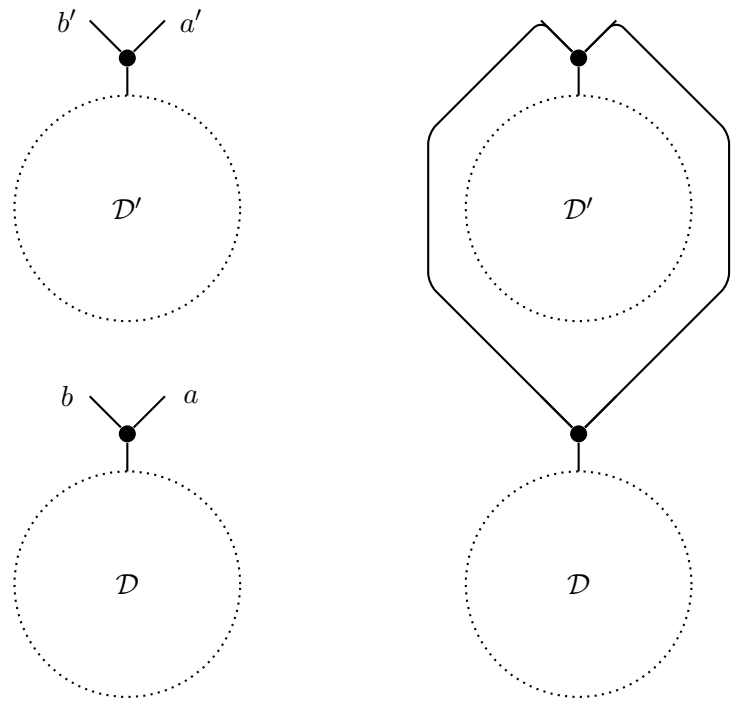
\begin{figure}[h!]
\begin{center}
\begin{tikzpicture}[scale=0.5, inner sep=0.8mm]

%%%%%%%

\draw [thick, dotted] (-3,10) arc (0:360:3);
\node (v2) at (-6,14) [shape=circle, fill=black] {};
\draw [thick] (-7,15) to (v2) to (-5,15);
\draw [thick] (v2) to (-6,13);

\node (a2) at (-4.3,15) {$a'$};
\node (b2) at (-7.6,15) {$b'$};
\node (D2) at (-6,10) {$\mathcal D'$};

%%%%%%

\draw [thick, dotted] (9,10) arc (0:360:3);
\node (v2') at (6,14) [shape=circle, fill=black] {};
\draw [thick] (7,15) to (v2') to (5,15);
\draw [thick] (v2') to (6,13);

\node (D2) at (6,10) {$\mathcal D'$};

%%%%%%%%%%

\draw [thick, dotted] (-3,0) arc (0:360:3);
\node (v1) at (-6,4) [shape=circle, fill=black] {};
\draw [thick] (-7,5) to (v1) to (-5,5);
\draw [thick] (v1) to (-6,3);

\node (a1) at (-4.3,5) {$a$};
\node (b1) at (-7.6,5) {$b$};
\node (D1) at (-6,0) {$\mathcal D$};

%%%%%%

\draw [thick, dotted] (9,0) arc (0:360:3);
\node (v1') at (6,4) [shape=circle, fill=black] {};
\draw [thick] (7,5) to (v1') to (5,5);
\draw [thick] (v1') to (6,3);

\node (D1) at (6,0) {$\mathcal D$};

%%%%%%

\draw [rounded corners, thick] (v1') to (10,8) to (10,12) to (7,15) to (v2');
\draw [rounded corners, thick] (v1') to (2,8) to (2,12) to (5,15) to (v2');

\end{tikzpicture}

\end{center}
\caption{$(1)$-composition of dessins}
\label{joining}
\end{figure}

If maps $\D$ and $\D'$ of degrees $n$ and $n'$ have $(k)$-handles $a, b$ and $a', b'$, for the same value of $k$, then $(a, a')$ and $(b, b')$ are defined to be new $2$-cycles for $x$, while the remaining cycles of $x$ on the two maps, together with those of $y$ and $t$, are unchanged. This process is called {\em $k$-composition}. Figure~\ref{joining} illustrates $(1)$-composition of $\D$ and $\D'$; the process is similar for $k=2, 3$. It is straightforward to check that $x$ and $y$ again generate a transitive group with $x^2=y^3=(xy)^7=1$, so this creates a new map of the required type and of degree $n+n'$. We will denote this new map by $\D+\D'$ or, following~\cite{Conder-80}, by $\D(k)\D'$ if we want to emphasise the type of handles used. If $v(\D)$ denotes the vector $(\alpha, \beta, \gamma)$ giving the numbers of fixed points of $x, y$ and $z$ in $\D$, then $v(\D+\D')=v(\D)+v(\D')-(4,0,0)$.

Topologically, this can be regarded as a connected sum operation on maps, in which two free edges $a, b$ of $\D$ are joined to two similar free edges $a', b'$ of $\D'$ to form two new edges, and the underlying surfaces are joined across cuts between the corresponding free ends. By drawing $\D$ and $\D'$ in the plane (possibly with crossings), both on the same axis of symmetry and with the relevant $(k)$-handles in the outer face, as on the left of Figure~\ref{joining}, one can see that no new crossings are introduced, so $\D+\D'$ has genus $g(\D+\D')=g(\D)+g(\D')$.

These composition operations can be iterated, provided suitable disjoint pairs of handles are available, allowing maps of unbounded degrees to be constructed. (Note that the only Hurwitz map containing non-disjoint handles is $\mathcal B$, where any two of its three handles have a fixed point in common.) The maps $\A,\ldots, \N$ corresponding to the basic diagrams in~\cite{Conder-80} all have genus $0$, so the same applies to any map constructed in this way by forming a chain (or even a tree) of them. However, by using other patterns of composition, including joining a map to itself, one can construct maps of unbounded genus: for instance the map $\G$, which has three $(1)$-handles, can act as a `pair of pants' (a sphere minus three discs), $2(g-1)$ of which can be composed with each other, using all available handles, to form an orientable surface of any genus $g\ge 2$. 

%%%%%%%%%%

\subsection{Cycles of $w$}

It will be important for the proof of Theorem~\ref{maintheorem} to have information about the cycles of the permutation $w:=xyt$ in various maps, so that Corollary~\ref{Jordancor} can be applied to the element $g=w$ and the subgroup $H:=\langle x, y\rangle$ of $G:=\langle x, y, t\rangle$.

The cycles of $w$ in the 14 basic diagrams are described in~\cite{Conder-80};  in addition to this one needs to understand the effect of $(k)$-composition on its cycles. For this it is sufficient to consider just those cycles containing $a, b, a'$ or $b'$, since no others are changed.

In $(1)$-composition, the two fixed points $a$ and $a'$ of $x$ in $\D$ and $\D'$, initially fixed by $w$ in its actions on $\D$ and $\D'$, form a $2$-cycle in its action on $\D(1)\D'$, while the cycles $(b, bw, \ldots, bw^{l-1})$ and $(b', b'w, \ldots, b'w^{l'-1})$ of $w$ containing $b$ and $b'$ merge to form a single cycle $(b, b'w, \ldots, b'w^{l'-1}, b', bw, \ldots, bw^{l-1})$ of $w$. More conveniently, if we write the two cycles as $c_b$ and $c_{b'}$ ending with $b$ and $b'$, they are concatenated to form a single cycle
\[c_b c_{b'}=(bw, \ldots, bw^{l-1}, b, b'w, \ldots, b'w^{l'-1}, b').\]

The situation is a little more complicated when $k=2$ or 3, depending on whether or not $a$ and $b$ are in the same cycle of $w$ on $\D$, and likewise for $a'$ and $b'$. If neither of these conditions is satisfied (as always happens when $k=1$), the cycles $c_a$ and $c_{a'}$ ending with $a$ and $a'$ are concatenated to form a cycle $c_ac_{a'}$, and similarly for the cycles $c_b$ and $c_{b'}$, as above. If $a$ and $b$ are in the same cycle of $w$, which we can write as $c_ac_b$, but $a'$ and $b'$ are not, then the cycles $c_{a'}$ and $c_{b'}$ are inserted, immediately after $a$ and $b$, into $c_ac_b$ to form a single cycle $c_ac_{a'}c_bc_{b'}$ of $w$. By symmetry between $\D$ and $\D'$ there is a similar insertion of cycles if $a'$ and $b'$ are in the same cycle, but $a$ and $b$ are not. However, if $a$ and $b$ are in the same cycle $c_ac_b$, and $a'$ and $b'$ are in the same cycle $c_{a'}c_{b'}$, then we obtain two cycles $c_ac_{b'}$ and $c_bc_{a'}$ of $w$ in $\D(k)\D'$. 

The following illustrative examples will be used later in the proof of Theorem~\ref{maintheorem}. First, for any map $\D$, let $\Pi(\D)$ denote the set of primes dividing the cycle-lengths of $w$ on $\D$.

\medskip

\noindent{\bf Example 1.} The map $\G$ has three $(1)$-handles, and the cycle structure of $w$ on $\G$ is $1^3, 13^3$ (see Table~\ref{BasicDessins}): the three fixed points are the elements $a$ in the three handles, and three cycles of length 13 each contain one of the points $b$ in the handles. Suppose that we compose two copies of $\G$ to form the map $\G(1)\G$, using a single $(1)$-handle $\{a, b\}$ and $\{a', b'\}$ in each copy. (It does not matter which handles we choose, by the $C_3$ symmetry of $\G$.) Then the two fixed points $a, a'$ form a $2$-cycle of $w$, and two cycles of length 13 are concatenated to form a $26$-cycle containing $b$ and $b'$. The cycle structure of $w$ on $\G(1)\G$ is therefore $1^4, 2, 13^4, 26$. This map now has four $(1)$-handles, and again each consists of a fixed point of $w$ and an element of a $13$-cycle. This means that if we iterate this process, adjoining further copies of $\G$ in any way to form a compound map $m\G=G+\cdots+\G$, the only cycle-lengths we will obtain will be $1, 2, 13$ and 26, so $\Pi(m\G)\subseteq\{2, 13\}$.

\medskip

\noindent{\bf Example 2.} The maps $\L$ and $\M$ each have a single $(2)$-handle, $\{a, b\}$ and $\{a', b'\}$. Suppose that we use these to form the map $\L(2)\M$. The cycle structure of $w$ on $\L$ is $21, 23, 58$, with $a$ and $b$ ($\alpha$ and $\beta$ in~\cite{Conder-80}) in the cycles of lengths 21 and 23. The cycle structure of $w$ on $\M$ is $1, 12, 14, 19, 26, 36$, with $a'$ and $b'$ in the cycles of lengths 36 and 19. Since $a$ and $b$ are in separate cycles of $w$, as are $a'$ and $b'$, the cycles containing $a$ and $a'$ are concatenated to form a cycle of length $21+36=57$, while those containing $b$ and $b'$ form a cycle of length $23+19=42$. Thus the cycle structure of $w$ on $\L(2)\M$ is $1, 12, 14, 26, 42, 57, 58$, so $\Pi(\L(2)\M)=\{2, 3, 7, 13, 19, 29\}$.

\medskip

In~\cite{Conder-80}, Conder defined a cycle $c$ of $w$ to be a {\em useful cycle\/} if, for each of the permutations $h=x$ and $y$, there is a point in $c$ with its image under $h$ also in $c$: this implies that $c$ satisfies condition~(2) of Corollary~\ref{Jordancor} with $H=\langle x, y\rangle\le G=\langle x, y, t\rangle$. It is also required that in the case $h=x$, the chosen point should not be a fixed point of $x$ contained in a $(k)$-handle for $k=1, 2$ or 3: this implies that the property of being useful is invariant under composition, even if this embeds the cycle in a longer cycle of $w$, since the latter is then also a useful cycle. In Table~\ref{BasicDessins}, the lengths of useful cycles in the 14 basic diagrams are indicated in bold face.

%%%%%%%%%%%%%%%

\section{Proof of Theorem~\ref{maintheorem}}\label{altgps}

The general strategy of the proof is, for each $r=0, 1, \ldots, 13$ and for each sufficiently large integer $m\equiv r$ mod~$(14)$, to construct a map $\W$ of degree $m$, by following Conder's method in~\cite{Conder-80}, and then to compose each of two copies of $\W$ with a map $\X_i$ ($i=1, 2$), where $\X_1$ and $\X_2$ have the same degree 210 but have different numbers of fixed points for each of the canonical generators $x, y$ and $z$. This gives two maps $\W_i=\W+\X_i$ of the same degree $n=m+210\equiv r$ mod~$(14)$, but with different numbers of fixed points for their canonical generators $x_i, y_i$ and $z_i$. As in~\cite{Conder-80} we will use Jordan's Theorem (or, more precisely, its corollary) to prove that the cartographic group $\langle x_i, y_i\rangle$ of each map $\W_i$ is the alternating group $A_n$. The condition on numbers of fixed points then guarantees that the two generating triples $(x_i, y_i, z_i)$ for $A_n$ satisfy the Beauville condition, with the Galois (minimal regular) covers $\widetilde{\W_i}$ of $\W_i$ as the corresponding Hurwitz curves $\C_i$ in the construction of the Beauville surface.

The first ingredient of  $\W$ is a map $\U_s$, corresponding to Conder's `stock' diagram $U_s$. Here we will make a slight modification of Conder's construction of $U_s$, for a reason which will be explained later. For each $s\ge 3$ we take $s'=\lfloor s/3\rfloor$ copies of the basic map $\G$, joined by $(1)$-handles to form a chain $s'\G=\G(1)\G(1)\cdots(1)\G$, and then if $s\equiv 1$ or $2$ mod $(3)$ we add a copy of the basic map $\A$ or $\E$ to this chain by $(1)$-composition. Since $\A$, $\E$ and $\G$ have degrees 14, 28 and 42, $\U_s$ has degree $14s$, and by inspection it has at least two unused $(1)$-handles, which we will need later. (This is why we have avoided the cases $s\le 2$ also considered by Conder, and why in the case $s\equiv 2$ mod 3 we have added a single copy of $\E$ rather than two copies of $\A$ as in~\cite{Conder-80}, thus ensuring that there are two unused $(1)$-handles when $s=5$.) As shown in Example~1, all cycles of $w$ in $s'\G$ have length $1, 2, 13$ or $26$, and in any unused $(1)$-handle, the fixed points $a$ and $b$ of $y$ are in cycles of $w$ of length $1$ and $13$ respectively. Adjoining $\A$ extends a pair of cycles of length 1 and 13 into cycles of length 2 and 26, while adjoining $\E$ increases their lengths to 2 and 22 and adds two cycles of length 9, so in all cases $\Pi(\U_s)\subseteq\{2, 3, 11, 13\}$.

The second ingredient of $\W$ is one of the maps corresponding to the fourteen compound diagrams listed by Conder in~\cite[p.~84]{Conder-80}; he did not name these diagrams (they are denoted by $H_d$ in~\cite{LTW}, where $d$ is their degree), so we shall call each map $\V_r$ for  $r=0, 1, \ldots, 13$ as it has degree $d_r\equiv r$ mod~$(14)$. These maps are constructed as chains of basic maps, described in the second column of Table~\ref{DessinsVr}. The cycle-lengths for $w$ on each map $\V_r$, calculated by Conder, are given in the final column, starting (before the semicolon) with the two cycles of lengths $1$ and $l_r$ containing the pair of points $a, b$ in an unused $(1)$-handle in $\V_r$. By using this $(1)$-handle, we can compose $\V_r$ with $\U_s$ for any $s\ge 3$ to form a map $\W=\U_s(1)\V_r$ of degree $m=14s+d_r\equiv r$ mod~$(14)$. The only effect of this composition on the cycle-lengths of $w$ is that two points $a$ in $\U_s$ and $a'$ in $\V_r$, originally fixed by $w$, now form a $2$-cycle, and two cycles of lengths $13$ and $l_r$ in $\U_s$ and $\V_r$, containing $b$ and $b'$, are joined to form a cycle of length $l'_r=13+l_r$, shown in the final column of Table~\ref{DessinsVr}. The lengths of the cycles of $w$ meeting $\V_r$ can therefore be seen from this table by replacing the initial pair $1, l_r$ with $2, l'_r$. The remaining cycles, those contained in $\U_s$, all have length 1, 2, 13 or 26 (and their multiplicities are unimportant). As shown by Conder, each map $\V_r$ provides a useful cycle of prime length $p_r$ (indicated in bold face in Table~\ref{DessinsVr}); by inspection of Table~\ref{DessinsVr}, $p_r$ is coprime to all other cycle-lengths for $w$ in $\W$.

\begin{table}[ht]
\centering
\begin{tabular}{| p{0.4cm} | p{3.5cm} | p{0.6cm} | p{5.9cm} | p{0.6cm} |}
\hline
$r$ & $\V_r$ & $d_r$ & Cycle-lengths for $w$ & $l'_r$ \\
\hline\hline
$0$ & $\H$ & $42$ & $1, 10;\, 3, 11, {\bf 17}$ & 23\\
\hline
$1$ & $\B(3)\H$ & $57$ & $1, 10;\, 3, {\bf 5}, 14, 24$ & 23 \\
\hline
$2$ & $\F(2)\E(1)\G(1)\H$& $142$ & $1, 13;\, 2^2, 3, 11, {\bf 17}, 22, 23, 24^2$ & 26 \\
\hline
$3$ & $\E(2)\I(2)\F$ & $115$ & $1, 9;\;\;\, 4, 10, {\bf 17}, 22^2, 30$ & 22 \\
\hline
$4$ & $\J(1)\K$ & $144$ & $1, 11;\, 2, 5, 10, 16, {\bf 17}, 22, 60$ & 24 \\
\hline
$5$ & $\C(3)\N(1)\E(2)\F$& $187$ & $1, 17;\, 2, 8, 18, 20, 24^2, 30, {\bf 43}$ & 30 \\
\hline
$6$ & $\B(3)\C(1)\G(1)\M(2)\F$ & 216 & $1, 13;\, 2^2, {\bf 5}, 8, 11, 12, 14, 24, 26, 34, 64$  & 26 \\
\hline
$7$ & $\C(1)\E(2)\E$ & $77$ & $1, 9;\;\;\, 2, 4, 8, {\bf 17}, 18^2$ & 22 \\
\hline
$8$ & $\B(3)\C$ & $36$ & $1, 11;\, {\bf 5}, 8, 11$ & 24 \\
\hline
$9$ & $\C(3)\H(1)\J$ & $135$ & $1, 11;\, 1, 2, 3, 8, 10, 16, {\bf 19}, 21^2, 22$ & 24 \\
\hline
$10$ & $\B(3)\C(1)\G(1)\E(2)\F$& $136$ & $1, 13;\, 2^2, {\bf 5}, 8, 11, 22, 24^3$ & 26 \\
\hline
$11$ & $\C(1)\J(1)\J$& $165$ & $1, 11;\, 2^2, 4, 8, 10^2, 16^2, {\bf 19}, 22^3$ & 24 \\
\hline
$12$ & $\J(1)\M$ & $180$ & $1, 11;\, 2, 10, 12, 14, 16, 19, 22, 26, {\bf 47}$ & 24 \\
\hline
$13$ & $\F(2)\I(2)\M$& $195$ & $1, 51;\, 4, 10, 12, 14, {\bf 23}, 26^2, 28$ & 64 \\
\hline

\end{tabular}
\caption{The maps $\V_r$}
\label{DessinsVr}
\end{table}

We will now modify these maps $\W$ by adjoining a suitable map $\X_i$ ($i=1,2$) to each of two copies of $\W$. We define $\X_1$ to be a chain $4\G=\G(1)\G(1)\G(1)\G$, with three copies of $\A$ attached by $(1)$-composition, so that $\X_1=4\G+3\A$ has degree $4\cdot 42 + 3\cdot 14 = 210\equiv 0$ mod~$(14)$, and has three unused $(1)$-handles. Since constructing $\X_1$ involves six joins, each eliminating four fixed points of $x$, the vector $v(\X_1)$ giving the numbers of fixed points of $x, y$ and $z$ on $\X_1$ is $4(6,0,0)+3(2,2,0)-6(4,0,0)=(6, 6, 0)$. Note that, as in the case of $\U_s$, all cycles of $w$ in $\X_1$ have length $1, 2, 13$ or $26$, so $\Pi(\X_1)=\{2, 3, 13\}$, and that in any unused $(1)$-handle, the fixed points $a$ and $b$ of $y$ are in cycles of $w$ of length $1$ and $13$ respectively. We define $\X_2$ to be the compound map $\L(2)\M$ of degree $102+108=210$ constructed in Example~2, with $\Pi(\X_2)=\{2, 3, 7, 13, 19, 29\}$. This has one unused $(1)$-handle, and its fixed point vector $v(\X_2)$ is $(2,0,4)+(4,0,3)-(4,0,0)=(2,0,7)$, differing from $v(\X_1)$ in all three coordinates. Since $\W$ has an unused $(1)$-handle (within $\U_s$), we can form a pair of maps $\W_i=\W(1)\X_i$ ($i=1,2$) of the same degree $n=m+210=14s+d_r+210\equiv r$ mod~$(14)$, but with fixed point vectors $v(\W_i)$ differing in all three coordinates. 

We now need to show that the cartographic group of $\W_i$ is the alternating group $A_n$. In order to apply Corollary~\ref{Jordancor}, we need the element $w$ to contain a useful cycle of prime length $p\le n-3$, which is coprime to all other cycle-lengths for $w$ in $\W_i$. The obvious candidate is the useful cycle of length $p=p_r$ in $\V_r$, since this is coprime to all other cycle-lengths of $w$ in $\W$.  We need to check that adjoining $\X_i$ does not introduce a new cycle which prevents this property from extending to $\W_i$.

Adjoining $\X_1$ to $\W$ creates new cycles of lengths 1, 2, 13 or 26, all of which are coprime to $p_r$ for each $r$. However, in the case of $\X_2$ we meet a technical problem: the cycle-lengths of $w$ on $\X_2=\L(2)\M$ are $1, 12, 14, 26, 23+19=42, 21+36=57$ and $58$, with the merged cycle of length $57$ containing the fixed point $b$ of $x$ in the unique $(1)$-handle; when we compose $\W$ and $\X_2$ by means of a $(1)$-handle in a copy of $\G$ in $\W_2$, we obtain a merged cycle of length $57+13=70$, which is divisible by 5. We also introduce prime divisors 3, 7, 19 and 29 from cycles in $\X_2$, which may or may not already be present in $\Pi(\W)$ (the primes 2 and 13 certainly are, from cycles in $\U_s$). The possibly new primes 3, 7 and 29 cause no problem since $p_r$ never takes these values. However, $p_r=5$ when $r=1, 6, 8$ or 10,  and $p_r=19$ when $r=9$ or 11, so in these six cases $p_r$ is not coprime to all other cycle-lengths of $w$ on $\W_2$.

%%%%%%%

\subsection{The cases $r=0, 2, 3, 4, 5, 7, 12, 13$.}

There is no such difficulty with the other eight values of $r$, those for which $p_r\in\{17, 23, 43, 47\}$. In these cases Corollary~\ref{Jordancor}, applied to the element $g=w$ and the transitive subgroup $H=\langle x, y\rangle$ of index $2$ in $G=\langle x, y, t\rangle$, proves that $G\ge A_n$, and hence $H=A_n$. For such values of $r$ this deals with all degrees $n=14s+d_r+210$ where $s\ge 3$.  For $r=0,  2, 3, 4, 5, 7, 12$ and $13$ we have 

\[d_r = 42,\, 142,\, 115,\, 144,\, 187,\, 77,\, 180\;\; {\rm and}\;\; 195,\]
so in these congruence classes $r$ mod~$(14)$ we have proved the result for all
\[n\ge 252+d_r = 294,\, 394, \,367,\, 396,\, 439,\, 329,\, 432\;\; {\rm and}\;\; 447.\]

%%%%%%%%

\subsection{The cases $r=6, 9, 10, 11$.}

When $r=6, 9, 10,$ or 11 the coprimality condition on $p_r$ fails for $j=2$, so the above argument does not apply. Let us therefore go back to the construction of $\W$ in these four cases, and replace the map $\V_r$ with $\V_{r^*}$ where $r^*=13, 2, 3, 4\equiv r+7$ mod~$(14)$, giving a map $\U_s(1)\V_{r^*}$, constructed as before, instead of $\U_s(1)\V_r$. This adds 7 to the congruence class mod~$(14)$ of the degree, so to restore the original class let us also adjoin a copy of $\C$ (of degree $21$) at an unused $(1)$-handle. If $s\ne 4, 5$ there are at least two such handles in $\U_s$, which can be used to attach $\C$ and $\X_j$, but if $s=4$ or $5$ there is only one, so we will need to create another by first adjoining an extra copy of $\G$ to $\U_s$. We therefore form $\W=\C(1)\U_{s^*}(1)\V_{r^*}$ and $\W_i=\W(1)\X_i$ for $i=1, 2$ where $s^*:=s+3$ or $s$ as $s=4,5$ or otherwise. This restores the congruence class of the degree $n$ mod~$(14)$, and replaces the prime $p_r=5$ with $p_{r^*}=23$ (when $r=6$) or 17 (when $r=9, 10$ or 11). There is a useful cycle of $w$ of length $p_{r^*}$ in $\W_i$; adjoining $\C$ creates merged cycles of length $1+1=2$ and $8+13=21$, and introduces new cycles of lengths 4 and 8, so $w$ has no other cycles of length divisible by $p_{r^*}$, and the proof proceeds as before. Since $n=21+14s^*+d_{r^*}+210$, with $s^*\ge 6$ and $d_{r^*}=195, 142, 115, 144$, we obtain lower bounds $n\ge 510, 457, 430$ and 459 for $r=6, 9, 10$ and 11.

%%%%%%%%%

\subsection{The cases $r=1, 8$.}

The trick we used for $r=6, 9, 10$ and 11, of replacing $\V_r$ with $\V_{r^*}$ where $r^*\equiv r+7$ mod~$(14)$, will not work in the remaining cases $r=1$ and 8, since it simply transposes these two problematic values, so here we must try something else. When $r=1$ let us instead replace $\V_1$ with $\V_5$ in the construction of $\W$, thus replacing $p_1=5$ with $p_5=43$. This adds 4 to the congruence class of the degree, so to restore it let us adjoin a copy of $\M$, of degree $108\equiv 10$ mod~$(14)$; this can be done using the $(1)$-handles in $\M$ and in a copy of $\G$ in $\U_{s^*}$ (where we define $s^*$ as before since we need an extra copy of $\G$ if $s=4$ or 5). Adjoining $\M$ to $\G$ creates cycles of $w$ of length $1+1=2$ and $36+13=49$, and introduces new cycles of length 12, 14, 19 and 26, but no cycles of length divisible by 43 arise. The proof now proceeds as before, giving a lower bound $n\ge 589$.

Similarly when $r=8$ we can replace $\V_8$ with $\V_{12}=\J(1)\M$, and $p_8=5$ with $p_{12}=47$. As in the case $r=1$, we can restore the congruence class of the degree by adding a copy of $\M$, only this time we can attach it, through its $(2)$-handle, to the unused $(2)$-handle in the copy of $\M$ in $\V_{12}$, so no extra copy of $\G$ is required. The useful cycle of $w$ of length $47=11+36$ in $\V_{12}$ is then merged with a cycle of length $36$ in $\M$ to form a cycle (still useful) of prime length $83$. There are no other cycles of length divisible by 83, so we can apply Corollary~\ref{Jordancor} as before. This gives a lower bound of $n\ge 540$, and completes the proof of Theorem~\ref{maintheorem}. \hfill$\square$

%%%%%%%%%%%%%%

\subsection{Smaller values of $n$}

Note that if $s=3$ we do not need a $(1)$-handle to adjoin $\A$ or $\E$ when constructing $\U_s$, so when $r=1, 6, 9, 10$ or $11$ we do not need to add an extra copy of $\G$ at the last stage  of the proof. In these cases we can therefore take $s^*=s=3$, reducing the degree by 42, so that $A_n$ is also a dHB group for $n=N_r-42=547$, $468$, $415$, $388$ and  $417$. (Note that we have not established the result for $n=N_r-14$ or $N_r-28$.)

For some values of $r$ we can find smaller examples of dHB groups by taking $s=0$ and using Conder's stock diagram $U_0$, namely the empty diagram, so that $\W=\V_r$. We then form $\W_i=\W(1)\X_i$ for $i=1,2$ by attaching ${\mathcal X}_i$ directly to $\V_r$, rather than to a copy of $\G$ in $\U_s$ as in the preceding proof, so that each $\W_i$ has degree $n=d_r+210$. (Here we use the $(1)$-handle in $\V_r$ which was previously used to attach it to $\U_s$.) To complete the proof, it is sufficient to check that the useful cycle of length $p_r$ in $\V_r$ remains useful after this joining operation. Attaching $\X_1$ merges the cycle of length $l_r$ in $\V_r$ containing a fixed point of $x$ with one of length $13$ in a copy of $\G$ in $\X_1$ to form a cycle of length $l_r'=l_r+13$; this is coprime to $p_r$ for each $r$, so no problem arises. However, attaching $\X_2$ creates a cycle of length $l_r''=l_r+57$, and now we have a problem when $r=4, 6$ or $10$, since then $l_r''=68, 70$ or $70$ is divisible by $p_r=17, 5$ or $5$. In all other cases, namely $r=0, 1, 2, 3, 5, 7, 8, 9, 11, 12$ or $13$, the prime $p_r$ does not divide $l_r''$, so $A_n$ is a dHB group for $n=252, 267, 352, 325, 397, 287, 246, 345, 375, 390$ and $405$.

In particular, $A_{246}$ is a dHB group, and in fact it is the smallest example of a dHB group which we can verify by hand. However, Marston Conder~\cite{Conder-17} has recently used Magma to show that $A_{168}$ is a dHB group. The following result shows that it is the smallest alternating group with this property. Indeed, it is currently the smallest example we know of any dHB group.

\begin{theorem}\label{168thm}
If $A_n$ is a dHB group then $n\ge 168$.
\end{theorem}

\noindent{\sl Proof.} If $A_n$ is a dHB group, then the Riemann-Hurwitz formula shows that the two quotient maps corresponding to its natural representation have genus $g_i$ satisfying
\[n=84(g_i-1) + 21\alpha_i+28\beta_i+36\gamma_i\quad(i=1,2),\]
where $\alpha_i, \beta_i$ and $\gamma_i$ are the numbers of fixed points of the two triples of generators $x_i, y_i$ and $z_i$. It follows that
\[n\equiv 21\alpha_i+28\beta_i+36\gamma_i\;{\rm mod}\, (84)\]
for $i=1,2$, so
\[\alpha_1\equiv\alpha_2\;{\rm mod}\,(4),\quad \beta_1\equiv\beta_2\;{\rm mod}\,(3),\quad \gamma_1\equiv\gamma_2\;{\rm mod}\,(7).\]
Writing
\[\alpha_1-\alpha_2=4a,\quad \beta_1-\beta_2=3b,\quad \gamma_1-\gamma_2=7c,\]
we see that $a, b$ and $c$ are integers satisfying
\[a+b+3c=g_2-g_1.\]
Moreover, they are each non-zero by the Beauville condition on the triples $(x_i, y_i, z_i)$. Since $g_i, \alpha_i, \beta_i, \gamma_i\ge 0$ it is clear that the least value of $n$ satisfying these conditions is obtained by taking each $g_i=0$, with $(a,b,c)=\pm(1,2,-1)$ or $\pm(2,1,-1)$, and taking the triples $(\alpha_i,\beta_i,\gamma_i)$ to be $(4,6,0)$ and $(0,0,7)$, or $(8,3,0)$ and $(0,0,7)$, so that $n=168$ in all cases. \hfill$\square$

%%%%%%%%%%%%%%%%%%%%%%%%%
 
\section{The double cover of $A_n$}\label{dblcov}

For each $n\ge 4$, the double cover $\tilde A_n=2.A_n$ is the unique group $G$ with a central involution $\iota\in G'$ such that $G/\langle \iota\rangle\cong A_n$. For $n=4$ or $5$ it is isomorphic to the binary tetrahedral or binary icosahedral group, and in the latter case also to $SL_2(5)$, while $\tilde A_6\cong SL_2(9)$.

An element of $A_n$ of odd order $k$ lifts to a pair of elements of order $k$ and $2k$ in $\tilde A_n$, while an involution $g\in A_n$ lifts to a pair of elements of order 2 or 4 as the number $\tau(g):=(n-|{\rm Fix}(g)|)/2$ of transpositions in $g$ (necessarily even) is or is not divisible by 4. It follows that a $(2,3,7)$-triple $(x,y,z)$ generating $A_n$ lifts to one generating $\tilde A_n$ if and only if $\tau(x)$ is divisible by 4, and that every such triple for $\tilde A_n$ arises in this way. Pellegrini and Tamburini~\cite{PT} used this to determine those $n$ for which $\tilde A_n$ is a Hurwitz group; these include all $n\ge 231$, together with some smaller values.

It also follows from a result of Fairbairn, Magaard and Parker~\cite{FMP} on quasisimple groups that $\tilde A_n$ is a Beauvillle group for all $n\ge 6$. Here we will show how to adapt the proof of Theorem~\ref{maintheorem} to prove Corollary~\ref{covercor}, that $\tilde A_n$ is a dHB group for all sufficiently large $n$.

It is sufficient to show that $A_n$ has a pair of generating triples $(x_i, y_i, z_i)$ of type $(2,3,7)$ for $i=1, 2$, which satisfy the Beauville property and have $\tau(x_i)$ divisible by 4, since these then lift to a suitable pair for $\tilde A_n$. We have already done all of this in proving Theorem~\ref{maintheorem}, apart from dealing with $\tau(x_i)$, so it is sufficient to adapt the proof to make each $\tau(x_i)/2$ even.

By inspection of Table~\ref{BasicDessins}, we see that in all 14 of Conder's basic diagrams, apart from $C, E, M$ and $N$, the involution $x$ has $\tau(x)/2$ odd. When forming a $(k)$-join, for any $k$, one must add the values of $\tau/2$ on the two components and also add 1 for the two new transpositions created by the join. Thus adjoining a copy of $\C$, $\E$, $\M$ or $\N$ to a map changes the parity of $\tau/2$, while adjoining any of the other ten basic maps preserves it. 

Since $\W_i=\W(1)\X_i$ for $i=1, 2$, with $\tau/2$ equal to $51$ and $52$ on $\X_1$ and $\X_2$, we see that $\tau(x_1)/2$ and $\tau(x_2)/2$ always have opposite parities. If $\tau(x_1)/2$ is odd and $\tau(x_2)/2$ is even, let us add a copy of $\E$ to $\W_1$, and two copies of $\A$ to $\W_2$, using $(1)$-joins with copies of $\G$ in $\U_s$; if the latter has insufficient unused $(1)$-handles we can add one or two copies of $\G$ to $\U_s$ in each $\W_i$. This changes the parity of $\tau(x_1)/2$ and preserves that of $\tau(x_2)/2$, so that both are now even, and it increases the degree $n$ in each case by 28 (or by 70 or 112  if extra copies of $\G$ are needed). It is straightforward to check that these adjunctions have no effect on the validity of the proof of Theorem~\ref{maintheorem}, other than possibly through changing the fixed point vectors $v(\W_i)$. Originally, we had $v(\W_1)-v(\W_2)=(6,6,0)-(2,0,7)=(4,6,-7)$, and adding copies of $\G$ to each $\W_i$ preserves this difference; adjoining $\E$ to $\W_1$ adds $(4,1,0)-(4,0,0)=(0,1,0)$ to $v(\W_1)$, and adjoining $2\A$ to $\W_2$ adds $2(2,2,0)-2(4,0,0)=(-4,4,0)$ to $v(\W_2)$, so after these adjunctions the  difference is $(8, 3, -7)$. This has all three coordinates non-zero, as required.

If, on the other hand, $\tau(x_1)/2$ is even and $\tau(x_2)/2$ is odd, we cannot simply transpose the roles of $\W_1$ and $\W_2$ since $v(\W_1)-v(\W_2)$ then becomes $(0,9,-7)$, so that $x_1$ and $x_2$ are conjugate in $A_n$. Instead, arranging as before that there are two unused $(1)$-handles in $\U_s$, we can use them to form a $(1)$-join in $\W_2$, thus increasing $\tau(x_2)/2$ by $1$, while leaving them unused in $\W_1$. This ensures that $\tau(x_i)/2$ is now even for $i=1, 2$. Again, the only significant effect on the proof of Theorem~\ref{maintheorem} is to increase the lower bounds on $n$ (now by 84), and to change the fixed point vectors $v(\W_i)$ of the two maps. Originally, we had $v(\W_1)-v(\W_2)=(6,6,0)-(2,0,7)=(4,6,-7)$; adding copies of $\G$ preserves this difference, and making an extra join in $\W_2$ subtracts $(4,0,0)$ from $v(\W_2)$, so the difference becomes $(8,6,-7)$, which is still non-zero in all three coordinates. This completes the proof of Corollary~\ref{covercor}.

%%%%%%%%%%%%%%%%%%%%%%%%%%%%%

\section{The groups $SL_n(q)$ and $L_n(q)$}\label{SLn}

In 2000, Lucchini, Tamburini and Wilson~\cite{LTW} showed that if $n\ge 287$ then $SL_n(q)$ is a Hurwitz group for each prime power $q$; it immediately follows that the same applies to the simple quotient group $L_n(q)$. Here we will combine their method of proof with ours for Theorem~\ref{maintheorem} to prove Corollary~\ref{SLnqcor}, that these groups are all dHB groups if $n$ is sufficiently large.

Following the notation of~\cite{LTW}, let $(\xi, y, (\xi y)^{-1})$ be the $(2,3,7)$-triple generating $A_n$ corresponding to the map $\W_i$ in the proof of Theorem~\ref{maintheorem}, where $n\ge N_r$ and $i=1$ or $2$. Adding an extra copy of $\G$ to $\U_s$ if necessary (and thus increasing the degree $n$ by 42), we may assume that the first copy of $\G$ in $\U_s$ has two free $(1)$-handles $a, b$ and $a', b'$ (these are $v_2, v_3$ and $v_{14}, v_{15}$ in~\cite{LTW}). For any prime power $q$ we can represent $\xi$ and $y$ as permutation matrices in $SL_n(q)$, acting naturally on the permutation module $V={\mathbb F}_q^n$. Let $t_1$ be a generator of the field $\mathbb F_q$, and let $x'$ (denoted by $x_1$ in~\cite{LTW}) be the matrix which acts on $a$ and $b$ by
\[a\mapsto -a+t_1a',\quad b\mapsto -b+t_1b',\]
while fixing the $n-2$ remaining standard basis vectors; thus $x'$ is an involution commuting with $\xi$. The results in~\cite[\S 4]{LTW}, applied there to Conder's generators $\xi, y$ for $A_n$, carry over to our situation without change, and show that if $x:=x'\xi$ and $z:=(xy)^{-1}$, then $(x,y,z)$ is a Hurwitz triple generating $SL_n(q)$.

Taking $i=1$ and 2 we obtain two such triples $(x_i, y_i, z_i)$. Now the subspace of $V$ fixed by the matrix $y_i$ has dimension equal to the number of cycles of the permutation $y_i$; this takes different values for $i=1$ and 2, since by our earlier construction the numbers of fixed points are different, so $y_1$ and $y_2$ cannot be conjugate in $SL_n(q)$ to powers of each other. A similar argument applies to the involutions $x_i$, since the relevant dimensions are the numbers of cycles of $\xi$, minus $2$. Finally, it is shown in~\cite[Lemma~4(a)]{LTW} that $xy$ is conjugate to $\xi y$, so this argument also applies to the elements $z_i$. Thus the triples $(x_i, y_i, z_i)$ satisfy the Beauville condition, as required. \hfill$\square$

\medskip

Lucchini and Tamburini~\cite{LT} have used the methods developed in~\cite{LTW} to prove Hurwitz generation for  
$Sp_{2n}(q)$, $\Omega^+_{2n}(q)$ and $SU_{2n}(q)$ for all $q$, and $\Omega_{2n+7}(q)$ and $SU_{2n+7}(q)$ for all odd $q$, provided $n\ge 371$ (and also for some smaller values of $n$). Similar arguments to that given above, omitted for brevity, show that these are also dHB groups for all sufficiently large $n$.

%%%%%%%%%%%%%%%%%%%%

\section{Other simple groups.}\label{simple}

It is natural to consider which other finite groups, especially simple groups, could be dHB groups. We will follow the ATLAS~\cite{Atlas-05} for notation for groups; in particular we will write $nA$, $nB$, $\ldots$ for the conjugacy classes of elements of order $n$ in a group, in decreasing order of their centralisers.

\begin{theorem}\label{sporadic}
No sporadic simple group is a doubly Hurwitz Beauville group.
\end{theorem}

\noindent{\sl Proof.} Twelve of the 26 sporadic simple groups are Hurwitz groups,  namely
\[J_1,\, J_2,\, J_4,\, Fi_{22},\, Fi_{24}',\, Co_3,\, He,\, Ru,\, HN,\, Ly,\, Th\;\;{\rm and}\;\; M\]
(see Conder~\cite{Conder-10}, for example); however, of these only Held's group $He$, Fischer's group $Fi_{24}'$ and the monster $M$ have at least two conjugacy classes of subgroups of each order $2$, $3$ and $7$, a necessary condition for being a dHB group. 

Held's group has two classes of involutions, two classes of elements of order 3 and five classes of elements of order 7. Its outer automorphism group, of order 2, transposes the mutually inverse classes $7A$ and $7B$, and similarly $7D$ and $7E$, whereas the classes $2A$, $2B$, $3A$, $3B$ and $7C$ are invariant; in looking for generating triples of type $(2,3,7)$ it is therefore sufficient to restrict attention to elements of order 7 in  classes $7A$, $7C$ or $7D$. As shown by Butler~\cite[\S3.4]{Butler-81}, the only such types with non-zero structure constants are $(2B,3A,7A)$, $(2A,3B,7C)$ and $(2B,3B,7D)$. Note that the Hurwitz triples for $He$  exhibited by Woldar in~\cite{Woldar-89} are of the last type. From Butler's analysis of the maximal subgroups of $He$ (specifically, his Theorem 23(i)), it follows by simple counting that triples of the first type all lie in maximal subgroups isomorphic to $S_4 \times L_3(2)$. This shows that any Hurwitz generating triple must have type $(2,3B,7)$, so $He$ cannot be a dHB group.

Linton and Wilson considered $Fi_{24}'$ in~\cite{LW-91}; they showed that generating $(2, 3, 7)$ triples must have type $(2A, 3C, 7A)$ or $(2A, 3E, 7B)$, which eliminates this group. Norton~\cite{Norton-98} showed that any generating $(2, 3, 7)$ triples in $M$ must be of type $(2B, 3B, 7B)$ or $(2B, 3C, 7B)$ (see also~\cite{Wilson-01}), so this group is also eliminated.  \hfill$\square$

\medskip

This leaves the simple groups of Lie type. Our aim now is to complete the proof of Theorem~\ref{simplethm} by eliminating various families of low rank. In~\cite[Prop.~3.1]{PT-IJGT}, Pellegrini and Tamburini have classified the simple Hurwitz groups with an absolutely irreducible projective representation of degree $n\le 7$. They are the following:

\begin{itemize}
\item $L_2(p)$ for primes $p\equiv 0, \pm 1$ mod~$(7)$ and $L_2(p^3)$ for $p\equiv\pm2, \pm 3$ mod~$(7)$;
\item $G_2(q)$ for prime powers $q\ge 5$;
\item `small' Ree groups ${}^2G_2(3^e)$ for odd $e>1$;
\item Janko's sporadic groups $J_1$ and $J_2$;
\item symplectic groups $S_6(q)=PSp_6(q)$ for odd $q\ge 5$; 
\item various groups $L_n(q)$ and unitary groups $PSU_n(q)$ for $n=5,6$ and $7$.
\end{itemize}
We have already shown that $J_1$ and $J_2$ are not dHB groups; in the following results we will also eliminate the first three of the above families, together with the `large' Ree groups ${}^2F_4(2^e)'$ for odd $e$ and the Steinberg groups ${}^3D_4(q)$ for odd $q$.

Many of these groups can immediately be eliminated because they do not have enough conjugacy classes of subgroups of order $2, 3$ or $7$:

\begin{lemma}\label{lowrank}
Let $G$ be a Hurwitz group belonging to one of the following families of finite simple groups:
\begin{enumerate}
\item the projective special linear groups $L_n(q)$ where $n\le 7$,
\item the small Ree groups $^2G_2(3^e)$ for odd $e > 1$,
\item the large Ree groups $^2F_4(2^e)'$ for odd $e \geq 1$,
\item the exceptional groups $G_2(q)$ for odd $q \geq 5$, or
\item the exceptional groups $^3D_4(q)$ for odd $q \geq 3$.
\end{enumerate}
Then $G$ is not a dHB group.\end{lemma}
\noindent{\sl Proof. } Many of these groups can be eliminated because they have a unique conjugacy class of involutions: this is clear for the groups $L_2(q)$, and has been proved by Ward~\cite{Ward-66} for the small Ree groups. For the groups $G_2(q)$ it has been proved by Chang~\cite[Theorem 4.4]{Chang-68} in characteristic $p\neq 2,3$, and by Enomoto~\cite{Enomoto-69} for $p=3$, while for $^3D_4(q)$ with $q$ odd the result is due to Deriziotis and Michler~\cite{DM-87}. The only Hurwitz group $L_n(q)$ with $n=3$ or $4$ is $L_3(2)\cong L_2(7)$, while Tamburini and Vsemirnov~\cite{TV-06} have shown that when $n=5, 6$ or $7$ the involutions in Hurwitz generating triples for $L_n(q)$ all lie in a single conjugacy class. Finally, Shinoda~\cite{Shinoda-75} has shown that the large Ree groups have a unique conjugacy class of elements of order 3.  \hfill$\square$

%%%%%%%%%%%%%%%%%%%%%%
	
\section{The groups $G_2(2^e)$}\label{G2}%G2(2^e)

In~\cite{Malle-1990}, Malle showed that the group $G_2(q)$ is a Hurwitz group if and only if $q\ge 5$. We dealt with odd values of $q$ in Lemma~\ref{lowrank}, so in this section we treat the remaining cases, where $q$ is even. Our aim is to prove the following:

\begin{theorem}\label{th:g2} If $G$ is the exceptional group $G_2(2^e)$ where $e \geq 3$, then $G$ is not a $dHB$ group.\end{theorem}

In order to prove the above theorem, it will be necessary to analyse the $7$-structure of $G_2(2^e)$, since the members of this family have two conjugacy classes each of elements of order 2 and 3. This result is due to Enomoto~\cite{Enomoto-69}, who also determined the structure of the centralisers of these elements, which we reproduce here.

\medskip

Let $q = 2^e$. Involutions in $G_2(q)$ belong to the conjugacy classes $2A$ or $2B$, with centralisers of order $q^6(q^2-1)$ or $q^4(q^2-1)$ respectively. The orders of the centralisers of elements of order $3$ depend on the congruence class $\epsilon = \pm 1$ of $q$ modulo $3$. We denote by $3A$ the class whose centralisers have order $q^3(q^2-1)(q^3-\epsilon)$ and are isomorphic to $SL_3(q)$ or $SU_3(q)$ respectively. 
We denote by $3B$ the class whose centralisers have order $q(q^2-1)(q-\epsilon)$ and are isomorphic to $L_2(q) \times (q-\epsilon)$. (Here $(q-\epsilon)$ is ATLAS notation for a cyclic group of order $q-\epsilon$.)

\medskip

The next lemma treats the conjugacy classes of cyclic subgroups of order~$7$ in $G_2(2^e)$.

\begin{lemma} \label{g2conj} Let $G = G_2(q)$ where $q = 2^e$ and $e \geq 1$.
\begin{enumerate}
\item If $q \equiv 2$ or $4$ {\rm mod}~$7$, or equivalently if $e$ is not divisible by $3$, then the Sylow $7$-subgroups of $G$ are cyclic.
\item If $q \equiv 1$ {\rm mod}~$7$, or equivalently if $e$ is divisible by $3$, then $G$ contains three conjugacy classes of cyclic subgroups of order~$7$.
\end{enumerate}
\end{lemma}

\noindent{\sl Proof. } The first part follows from the list of conjugacy classes of $G$ in~\cite{Enomoto-69}. The second follows again from \cite{Enomoto-69} and from Malle's analysis of the Sylow $7$-subgroups of $G$ in~\cite{Malle-1990}.\hfill$\square$

\medskip

As an immediate corollary we have the following.

\begin{corollary} If $2^e\equiv 2$ or $4$ \rm mod~$7$, then $G_2(2^e)$ is not a $dHB$ group.\end{corollary}

\medskip

We may therefore assume from now on that $e$ is divisible by $3$, so we will change notation and write $q=8^e$. Note that $\epsilon = 1$ or $-1$ as $e$ is even or odd.

Following Malle's analysis~\cite{Malle-1990} we see that one class of cyclic subgroups of order~$7$ in $G$ are generated by elements from a unique conjugacy class of regular elements which we denote by $7C$. The remaining two classes of cyclic subgroups of order~$7$ are generated by elements from one of two conjugacy classes which we denote by $7A$ and $7B$. Elements from these two classes are distinguished by the fact that if $s \in 7A$ or $7B$ and $t$ is an involution in $C_G(s)$, then $t \in 2A$ or $2B$ respectively~\cite{Enomoto-69}. If $s \in 7A \cup 7B$ then $C_G(s)\cong SL_2(q) \times (q-1)$, whereas $C_G(s)\cong (q-1)\times(q-1)$ if $s\in7C$.

\begin{remark} The Hurwitz triples given by Malle \cite{Malle-1990} are of type $(2B,3B,7C)$ in our notation.\end{remark}

\medskip

Our strategy is to show that there do not exist Hurwitz generating triples for $G_2(8^e)$ of type $(2A,3,7)$ by using the structure constants for $G$, as described in Theorem~\ref{th:frobenius}. We begin by determining the non-zero $(2A,3,7)$ structure constants in CHEVIE~\cite{CHEVIE}, which are as follows:
\[n(2A,3B,7B) = q^7(q^6-1), \quad n(2A,3B,7C) = \vert G \vert\]

We immediately have the following:

\begin{lemma} \label{abb} Let $G = G_2(q)$ where $q=8^e$ and let $H$ be a Hurwitz subgroup of $G$ of type $(2A,3B,7B)$. Then $H$ is a proper subgroup of $G$.\end{lemma}

\noindent{\sl Proof. } Define the following set
\[\H = \{ (x,y,z) \in G \times G \times G \mid x \in 2A,\, y \in 3B,\, z \in 7B,\, xyz=1\} \]
and note that $\vert \H \vert = n(2A,3B,7B) = q^7(q^6-1)$. The Hurwitz subgroup $H$ is generated by a triple $(x,y,z) \in \H$. We now let each $g \in G$ act on $\H$ by sending $(x,y,z)$ to $(x^g,y^g,z^g)$. Since $\vert G \vert = q^6(q^6-1)(q^2-1) > \vert \H \vert $, any orbit of $G$ on $\H$ has non-trivial stabilisers, that is, $H$ has a non-trivial centraliser. It follows, since $G$ is non-abelian and simple, that $H$ is a proper subgroup of $G$.\hfill$\square$

\medskip

For the Hurwitz subgroups of type $(2A,3B,7C)$, we must be a little more delicate.

\begin{lemma} \label{abc} Let $G = G_2(q)$ where $q=8^e$, and let $H$ be a Hurwitz subgroup of $G$ of type $(2A,3B,7C)$.\begin{enumerate}
\item If $\epsilon = -1$ then $H$ lies in a unique conjugacy class of subgroups $H \cong L_3(2)$, with $N_G(H) \cong L_3(2).2$.
\item If $\epsilon = 1$ then $H$ lies in one of two conjugacy classes of subgroups isomorphic to $L_3(2)$, with representatives $H_1$ and $H_2$ where $N_G(H_1) \cong 3 \times L_3(2)$ and $N_G(H_2) \cong (3 \times L_3(2)).2$.
\end{enumerate}
\end{lemma}

\noindent{\sl Proof. } First we show that such subgroups $H$ exist. From the work of Cooperstein~\cite{Cooperstein-1981} we see that there exists a unique conjugacy class of maximal subgroups of $G$ isomorphic to $SL_3(q).2$. The commutator subgroup $SL_3(q)$ of such a maximal subgroup contains a unique conjugacy class of involutions whose elements belong to the $G$-class $2A$. The classification of maximal subgroups of $SL_3(q)$ can be found in~\cite{BHR-2013}, from which we see that such groups contain subgroups isomorphic to $L_3(2)$. It is clear that elements of order $7$ in such subgroups belong to the $G$-class $7C$, and from the non-zero structure constants of $G$, elements of order $3$ must then belong to the $G$-class $3B$.

From the list of maximal subgroups of $SL_3(q)$ again we see that there are $d=$ gcd$(q-1,3)$ conjugacy classes of subgroups of shape $d \times L_3(2)$ in $SL_3(q)$. Suppose that $d=1$, in which case $\epsilon = -1$, and let $H \cong L_3(2)$ be a representative of such a conjugacy class. The number of $(2,3,7)$ triples in $H$ can be counted by using Theorem~\ref{th:frobenius}; there are $2 \vert H \vert$ of them, since we must count across both conjugacy classes of elements of order $7$ in $H$. The normaliser in $G$ of $H$ is $N_G(H) \cong L_3(2).2$, so $H$ has $|G:N_G(H)|=|G|/2|H|$ conjugates, giving $|G|$ triples. this accounts for the full number $n_G(2A,3B,7C)$ of triples. Now suppose that $d=3$, in which case $\epsilon = 1$. There are three conjugacy classes of subgroups isomorphic to $L_3(2)$ in $SL_3(q)$, two of which are fused under the action of the outer graph automorphism. Denote by $H_1$ and $H_2$ conjugacy representatives of these classes such that $N_G(H_1) \cong L_3(2) \times 3$ and $N_G(H_2) \cong (3 \times L_3(2)).2$. A similar calculation shows that the contribution from these classes is equal to $\vert G \vert$. This completes the proof.\hfill$\square$

\medskip

We have shown that any Hurwitz generating triple for $G$ is of type $(2B,3,7)$, completing the proof of Theorem \ref{th:g2}.

%%%%%%%%%%%%%%%%%%%%%%

\section{The groups $^3D_4(2^e)$}\label{3D4}

%3D4(2^e)

The groups $^3D_4(q)$ contain natural copies of $G_2(q)$ as subgroups, so many of the results we need for them, together with the techniques we use, are similar to those for $G_2(q)$. Malle~\cite{Malle-1995} has shown that $^3D_4(p^e)$ ($p$ prime) is a Hurwitz group if and only if $p \neq 3$ and $p^e \neq 4$. The case where $p \geq 5$ has already been dealt with in Theorem~\ref{lowrank}, so our aim in this section is to prove the following:

\begin{theorem} \label{th:3d4} If $G$ is the Steinberg triality group $^3D_4(2^e)$ where $e \geq 1$, then $G$ is not a dHB gorup.\end{theorem}

As in the case of the groups of type $G_2(2^e)$, the groups of type $^3D_4(2^e)$ have two conjugacy classes each of elements of order $2$ and $3$, as determined by Deriziotis and Michler \cite{DM-87}. Here we reproduce their results on the centralisers of these elements.

Let $q=2^e$. Involutions in $^3D_4(q)$ belong to the conjugacy classes $2A$ or $2B$, with centralisers of order $q^{12}(q^6-1)$ or $q^{10}(q^2-1)$ respectively. The centralisers of elements of order $3$ again depend on the congruence class $\epsilon = \pm1$ of $q$ modulo $3$. We denote by $3A$ the class whose elements have centralisers isomorphic to $SL_2(q^3) \times (q-\epsilon)$, where $(q - \epsilon)$ denotes the cyclic group of that order, as before. We denote by $3B$ the class for which the centralisers are isomorphic to $GL_3(q)$ or $GU_3(q)$ when $\epsilon = 1$ or $-1$ respectively.

We now turn to the $7$-structure of the groups $^3D_4(2^e)$.

\begin{lemma} Let $G = {}^3D_4(q)$ where $q=2^e$ and $e \geq 1$.\begin{enumerate}
\item If $q \equiv 2$ or $4$ {\rm mod} $7$, or equivalently if $e$ is not divisible by $3$, then $G$ has two conjugacy classes of cyclic subgroups of order $7$.
\item If $q \equiv 1$ {\rm mod} $7$, or equivalently if $e$ is divisible by $3$, then $G$ has three conjugacy classes of cyclic subgroups of order $7$.
\end{enumerate}\end{lemma}
\noindent{\sl Proof. } These results follow from the analysis due to Malle \cite{Malle-1995} and the determination of the conjugacy classes of $^3D_4(q)$ in \cite{DM-87}.\hfill$\square$

\medskip

We are now in a position to prove the following.

\begin{lemma} If $2^e \equiv 2$ or $4$ {\rm mod} $7$, then $^3D_4(2^e)$ is not a dHB group.\end{lemma}
\noindent{\sl Proof. } The structure constants for $^3D_4(2^e)$ in these cases can be determined by using CHEVIE. We see that all of the $n(2A,3,7)$ structure constants are zero, hence there cannot exist a compatible dHB structure for $G$.\hfill$\square$

\medskip

It now suffices to consider the case where $e$ is divisible by $3$, so we change notation and write $q=8^e$. Note that $\epsilon=1$ or $-1$ as $\epsilon$ is even or odd. We now provide a more detailed description of the $7$-structure of $^3D_4(8^e)$ following \cite{DM-87} and Malle's analysis \cite{Malle-1995}.

There exists a unique conjugacy class of regular elements of order $7$ in $^3D_4(q)$, which we denote by $7C$. If $s \in 7C$, then $C_G(s) \cong (q^3-1) \times (q-1)$. The remaining two classes of cyclic subgroups of order 7 are generated by elements from one of two conjugacy classes, which we denote by 7A and 7B. If $s \in 7A$ then $C_G(s) \cong SL_2(q^3) \times (q-1)$, while if $s \in 7B$ then $C_G(s) \cong SL_2(q) \times (q^3-1)$.

\begin{remark} In our notation, the Hurwitz triples shown to exist by Malle~\cite{Malle-1995} are of type $(2B,3B,7C)$.\end{remark}

We are now in a position to prove the following lemma, which will complete the proof of Theorem \ref{th:3d4}.

\begin{lemma} \label{aac} Let $G = {}^3D_4(q)$ where $q=8^e$, and let $H$ be a Hurwitz subgroup of $G$ of type $(2A,3,7)$. Then $H$ is a proper subgroup of $G$.\end{lemma}
\noindent{\sl Proof. } We must first determine the non-zero values of $n(2A,3,7)$ for $G$; there are two such candidates: $n_G(2A,3A,7A)$ and $n_G(2A,3A,7C)$.

Using CHEVIE we calculate
\[n_G(2A,3A,7A) = \frac{q^3-q^2+1}{q^2(q^6-1)}\vert G \vert.\]
By an argument similar to that used in the proof of Lemma \ref{abb}, it follows that Hurwitz subgroups of $G$ of type $(2A,3A,7A)$ are always proper.

The value $n(2A,3A,7C)$ is a so-called `exceptional value' in CHEVIE, and so it (incorrectly) evaluates to 0 despite contributions from subgroups isomorphic to $G_2(q)$. Hence we must determine the value by hand. For this we use the character table of $G$ to determine a bound on the order of $n(2A,3A,7C)$. The unipotent characters were determined by Spaltenstein~\cite{Spaltenstein-82} and the remaining characters were determined by Deriziotis and Michler \cite{DM-87}. For simplicity, we assume that $\epsilon =-1$; the proof in the case $\epsilon = 1$ is identical, apart from appropriate sign changes in the character values. Apart from the trivial character, the characters which do not take the value $0$ on all three conjugacy classes are those denoted as follows: $[\varepsilon_1]$, $[\varepsilon_2]$ and the four families of characters $\chi_{_{4,1}}$, $\chi_{_{4,St}}$, $\chi_{_{5,1}}$ and $\chi_{_{5,St}}$. Crucially, since any element $g \in 2A \cup 3A \cup 7C$ is conjugate to all of its non-trivial powers, the character values of $g$ are all rational -- in particular they are all integers -- which allows us to compute the following bound relatively easily:
\[n_G(2A,3A,7C) \leq \frac{|G|}{q^2+q+1}.\]
From this we see that any $(2A,3A,7C)$-subgroup has a non-trivial centraliser in $G$, completing the proof.\hfill$\square$

%%%%%%%%%%%%%%%%%%%%%

 \section{Connections with maps and dessins d'enfants}\label{dessins}

In this section we will put the specific arguments and constructions used in Section~\ref{altgps}, and those in~\cite{Conder-80} on which they are based, in the more general context of the algebraic theory of maps on surfaces~\cite{JS-78} and of Grothendieck's theory of dessins d'enfants~\cite{Gro}. For further background, see~\cite{JW-16}.

Let $G$ be a Hurwitz group, generated by a Hurwitz triple $(x,y,z)$, and let $K=\ker\theta$ where $\theta$ is the epimorphism $\Delta\to G$ given by $X\mapsto x, Y\mapsto y$ and $Z\mapsto z$. There is a natural action of $\Delta$ as a group of orientation-preserving isometries of the hyperbolic plane $\mathbb H$, or equivalently as a group of automorphisms of the Riemann surface $\mathbb H$;  this action preserves the universal map $\U$ of type $\{7, 3\}$, a trivalent tessellation of $\mathbb H$ by hyperbolic heptagons. This induces an action of $\Delta/K$, and hence of $G$, as the automorphism group of the compact Riemann surface $\R=\mathbb H/K$ of genus $\widetilde g=1+|G|/84$, and also as the orientation-preserving automorphism group of $\U/K$, an orientably regular map of type $\{7,3\}$ on $\R$.
 
If $G$ has a faithful, transitive permutation representation of degree $n$, with $H$ the stabiliser of a point, let $M=\theta^{-1}(H)$, a subgroup of index $n$ in $\Delta$ with core $K$. Then the Riemann surface $\mathcal S=\mathbb H/M=\R/H$ has a tessellation $\T=\U/M$, a map with $\U/K$ as its Galois (minimal regular) cover $\widetilde\T$ and $G$ as its monodromy group, that is, the monodromy group of the branched covering $\R\to\mathcal S$.
 
Having finite index in $\Delta$, both $K$ and $M$ are finitely generated Fuchsian groups, with $K$ isomorphic to the fundamental group of $\R$, a surface group of genus $\widetilde g$. By a theorem of Singerman~\cite{Singerman-70} {(see also \S\ref{classical}), $M$ has signature
\[(g; 2^{[\alpha]}, 3^{[\beta]}, 7^{[\gamma]})\]
where there are $\alpha=|{\rm Fix}(x)|$ elliptic periods equal to $2$, $\beta=|{\rm Fix}(y)|$ equal to $3$, and  $\gamma=|{\rm Fix}(z)|$ equal to $7$, and where the genus $g$ (of $M$, $\mathcal S$ and $\T$) is given by applying the Riemann-Hurwitz formula~(\ref{RHeqn}) to the inclusion $M\le\Delta$.

\smallskip

Being compact, the Riemann surfaces $\R$, $\mathcal S$ and $\mathbb H/\Delta$ can be regarded as complex algebraic curves. The latter, having genus $0$, can be identified with the Riemann sphere, or complex projective line, $\hat\mathbb C=\mathbb P^1(\mathbb C)=\mathbb C\cup\{\infty\}$. The inclusions of $M$ and $K$ in $\Delta$ induce meromorphic functions $\beta:\mathcal S\to\hat\mathbb C$ and $\widetilde\beta:\R\to\hat\mathbb C$ which are branched over three points (the images of the fixed points of $X, Y$ and $Z$, which can be identified with $0, 1$ and $\infty$), so by Bely\u\i's Theorem~\cite{Belyi-79} the pairs $(\mathcal S,\beta)$ and $(\R,\widetilde\beta)$ are defined over the field $\overline\Q$ of algebraic numbers. Indeed, it follows from recent work of Gonz\'alez-Diez and Jaikin-Zapirain~\cite{GDJZ} (see also Kucharczyk~\cite{Kucharczyk-14}) that these regular Bely\u\i\/ pairs $(\R,\widetilde\beta)$ associated with Hurwitz groups provide a faithful representation of the absolute Galois group ${\rm Gal}\,\overline\Q/\Q$.

The maps $\T$ and $\widetilde\T$ are the dessins d'enfants associated with the Bely\u\i\/ pairs $(\mathcal S,\beta)$ and $(\R,\widetilde\beta)$, that is, they are the inverse images under $\beta$ and $\widetilde\beta$ of the unit interval $[0, 1]\subset\hat\mathbb C$, with the vertices forming the fibre over $0$. If $\T$ is one of the basic maps $\A,\ldots,\N$, or one of those constructed from them in Section~\ref{altgps} or (implicitly) in~\cite{Conder-80}, then $\mathcal S$ has genus $0$, so $\beta$ is a rational function $\hat\mathbb C\to\hat\mathbb C$ of degree $n$. However, explicitly determining $\beta$, even in a simple case such as $\mathcal S=\A$, where $n=14$ and $G\cong L_2(13)$, is a difficult problem.

\medskip

\noindent{\bf Example} Macbeath~\cite{Macbeath-69} showed that for each prime $p\equiv\pm 1$ mod~$(7)$ there are three normal subgroups $K_i\;(i=1, 2, 3)$ of $\Delta$ with $\Delta/K_i\cong L_2(p)$, corresponding to three Riemann surfaces $\R_i$ of genus $\widetilde g=1+p(p^2-1)/168$. In the corresponding Hurwitz triples generating $G$, the elements $x$ and $y$ lie in the unique conjugacy classes of elements of order 2 and 3, while the elements $z$ lie in the three different conjugacy classes of elements of order 7. As shown by Streit~\cite{Streit-00}, the corresponding Bely\u\i\ pairs $(\R_i,\widetilde\beta)$ are defined over the real cyclotomic field $\Q(\cos(2\pi/7))$, and are conjugate under its Galois group $C_3$. If we take the natural representation of $G$ of degree $n=p+1$, with $H\cong C_p\rtimes C_{(p-1)/2}$ a Sylow $p$-normaliser, we obtain three quotient surfaces $\mathcal S_i=\R_i/H$. Their common genus $g$ can be found by using the fact that the numbers $\alpha, \beta$ and $\gamma$ of fixed points of $x, y$ and $z$ are 0 or 2 as their orders 2, 3 and 7 divide $(p\pm 1)/2$ respectively. For instance, if $p=13$ we have $\alpha=\beta=2$ and $\gamma=0$, so $g=0$; one of the three planar maps $\mathcal S_i$ is the map $\A$ shown in Figure~\ref{DiagramA}. For $p=29, 41, 43,\ldots$ we have $g=0, 1, 0,\ldots\sim p/84$. 

\medskip
 
The realisations of various groups $G=A_n$ as Hurwitz groups in Section~\ref{altgps} (and also, implicitly, in~\cite{Conder-80}) are all constructed in the above way, using the natural representation of $G$ of degree $n$, with permutations $x, y$ and $z$ defining maps $\W_i$ which play the role of $\T$ and are designed to have monodromy group $\langle x, y, z\rangle\cong A_n$. In each of these examples, $\T$ has genus $0$ (although not strictly necessary, this simplifies exposition and illustration), so that $n=21\alpha+28\beta+36\gamma-84$. By contrast, in the last paragraph of Section~\ref{dblcov}, the extra $(1)$-join created in $\W_2$ increases its genus from $0$ to $1$ while reducing $\alpha$ by $4$, so that $n=21\alpha+28\beta+36\gamma$.

Similarly, for the construction of Hurwitz groups it is convenient though not really necessary that the maps $\T$ and therefore also $\widetilde\T$, both here and in~\cite{Conder-80}, are invariant under a reflection $t$. (Conder's motivation for this was also to realise various symmetric and alternating groups as quotients of the extended triangle group $\Delta^*=\Delta[2,3,7]$.) It follows that in such cases $K$ is normal in $\Delta^*$, and the full automorphism group of the map $\widetilde\T$, including orientation-reversing elements, is $\Delta^*/K$. When $G=A_n$ this is isomorphic to either $A_n\times C_2$ or $S_n$ as the permutation $t$ is even or odd; in the examples in~\cite{Conder-80} and Section~\ref{altgps}, the parity of $t$ can be determined from the information in the third column of Table~\ref{BasicDessins}.

 %%%%%%%%%%%%%%%%%%%%%%%%%

%%%%%%%%%%%%%%%%%%%%%

\newpage

\section{Appendix: The 14 basic maps}

Here we show the maps $\A,\ldots, \N$ corresponding to Conder's basic diagrams $A\ldots, N$.

\bigskip

\begin{figure}[h]
\begin{center}
\begin{tikzpicture}[scale=0.5, inner sep=0.8mm]

\draw [thick] (2,0) arc (0:360:2);
\node (a) at (2,0) [shape=circle, fill=black] {};
\node (a') at (-2,0) [shape=circle, fill=black] {};
\node (b) at (0,-2) [shape=circle, fill=black] {};
\node (c) at (0,-0.5) [shape=circle, fill=black] {};
\node (d) at (4.5,0) [shape=circle, fill=black] {};
\node (d') at (-4.5,0) [shape=circle, fill=black] {};

\draw [thick] (b) to (c);
\draw [thick] (1,0.5) to (c) to (-1,0.5);
\draw [thick] (a) to (d);
\draw [thick] (a') to (d');

\node at (0,-4)  {$\A$};

%%%%%%

\draw [thick] (12,0) arc (0:360:2);
\node (b2) at (10,-2) [shape=circle, fill=black] {};
\node (c2) at (11,-1.7) [shape=circle, fill=black] {};
\node (c'2) at (9,-1.7) [shape=circle, fill=black] {};
\node (d2) at (10,2) [shape=circle, fill=black] {};

\draw [thick] (c2) to (11,0);
\draw [thick] (c'2) to (9,0);
\draw [thick] (d2) to (10,0);

\node (e2) at (10,-4) [shape=circle, fill=black] {};
\draw [thick] (b2) to (e2);
\draw [thick] (12,0) arc (0:360:2);
\draw [thick] (10.5,-4.5) arc (0:360:0.5);

\node at (12,-4)  {$\B$};

\end{tikzpicture}

\end{center}
\caption{The maps $\A$ and $\B$.}
\label{dessinsAB}
\end{figure}

%%%%%%%%%%%%%%%%

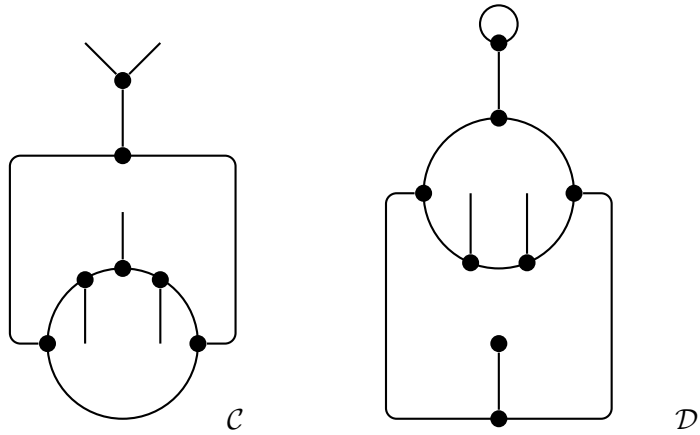
\begin{figure}[h!]
\begin{center}
\begin{tikzpicture}[scale=0.5, inner sep=0.8mm]

\draw [thick] (2,0) arc (0:360:2);
\node (a) at (2,0) [shape=circle, fill=black] {};
\node (a') at (-2,0) [shape=circle, fill=black] {};
\node (b) at (0,2) [shape=circle, fill=black] {};
\node (c) at (1,1.7) [shape=circle, fill=black] {};
\node (c') at (-1,1.7) [shape=circle, fill=black] {};

\draw [thick] (c) to (1,0);
\draw [thick] (c') to (-1,0);
\draw [thick] (b) to (0,3.5);

\node (d) at (0,5) [shape=circle, fill=black] {};
\draw [rounded corners, thick] (a) to (3,0) to (3,5) to (-3,5) to (-3,0) to (a');

\node (e) at (0,7) [shape=circle, fill=black] {};
\draw [thick] (d) to (e);
\draw [thick] (1,8) to (e) to (-1,8);

\node at (3,-2)  {$\C$};

%%%%%%%

\draw [thick] (12,4) arc (0:360:2);
\node (a4) at (12,4) [shape=circle, fill=black] {};
\node (a'4) at (8,4) [shape=circle, fill=black] {};
\node (b4) at (10,6) [shape=circle, fill=black] {};
\node (c4) at (10.75,2.15) [shape=circle, fill=black] {};
\node (c'4) at (9.25,2.15) [shape=circle, fill=black] {};
\node (d4) at (10,8) [shape=circle, fill=black] {};
\node (e4) at (10,0) [shape=circle, fill=black] {};
\node (f4) at (10,-2) [shape=circle, fill=black] {};

\draw [thick] (c4) to (10.75,4);
\draw [thick] (c'4) to (9.25,4);
\draw [thick] (b4) to (d4);
\draw [thick] (e4) to (f4);
\draw [rounded corners, thick] (a4) to (13,4) to (13,-2) to (7,-2) to (7,4) to (a'4);

\draw [thick] (10.5,8.5) arc (0:360:0.5);

\node at (15,-2)  {$\D$};

\end{tikzpicture}

\end{center}
\caption{The maps $\C$ and $\D$.}
\label{dessinsCD}
\end{figure}

\newpage

%%%%%%%%%%%%%%%%%%%%

\begin{figure}[h!]
\begin{center}
\begin{tikzpicture}[scale=0.45, inner sep=0.8mm]

\draw [thick] (2,5) arc (0:360:2);
\draw [thick] (2,0) arc (0:360:2);
\node (a) at (2,5) [shape=circle, fill=black] {};
\node (a') at (-2,5) [shape=circle, fill=black] {};
\node (b) at (0,3) [shape=circle, fill=black] {};
\node (c) at (0,4.5) [shape=circle, fill=black] {};
\node (d) at (0.8,1.8) [shape=circle, fill=black] {};
\node (d') at (-0.8,1.8) [shape=circle, fill=black] {};
\node (e) at (2,0) [shape=circle, fill=black] {};
\node (e') at (-2,0) [shape=circle, fill=black] {};
\node (f) at (0,-2) [shape=circle, fill=black] {};
\node (g) at (0,-3.5) [shape=circle, fill=black] {};

\draw [rounded corners, thick] (a) to (3,5) to (3,0) to (e);
\draw [rounded corners, thick] (a') to (-3,5) to (-3,0) to (e');
\draw [thick] (b) to (c);
\draw [thick] (1,5.5) to (c) to (-1,5.5);
\draw [thick] (d) to (0.8,0);
\draw [thick] (d') to (-0.8,0);
\draw [thick] (f) to (g);

\node at (3,-3.5) {$\E$};

%%%%%%

\draw [thick] (17,0) arc (0:360:2);
\node (a6) at (17,0) [shape=circle, fill=black] {};
\node (a'6) at (13,0) [shape=circle, fill=black] {};
\node (b6) at (15,2) [shape=circle, fill=black] {};
\node (c6) at (15.75,-1.85) [shape=circle, fill=black] {};
\node (c'6) at (14.25,-1.85) [shape=circle, fill=black] {};

\draw [thick] (c6) to (15.75,0);
\draw [thick] (c'6) to (14.25,0);

\node (d6) at (15,4) [shape=circle, fill=black] {};
\draw [thick] (b6) to (d6);
\node (e6) at (19,4) [shape=circle, fill=black] {};
\node (e'6) at (11,4) [shape=circle, fill=black] {};
\node (f6) at (19,2) [shape=circle, fill=black] {};
\node (f'6) at (11,2) [shape=circle, fill=black] {};

\draw [thick] (19.5,1.5) arc (0:360:0.5);
\draw [thick] (11.5,1.5) arc (0:360:0.5);

\draw [rounded corners, thick] (a6) to (21,0) to (21,4) to (9,4) to (9,0) to (a'6);
\draw [thick] (e6) to (f6);
\draw [thick] (e'6) to (f'6);

\node at (18,-3.5) {$\F$};

\end{tikzpicture}

\end{center}
\caption{The maps $\E$ and $\F$.}
\label{DessinsEF}
\end{figure}
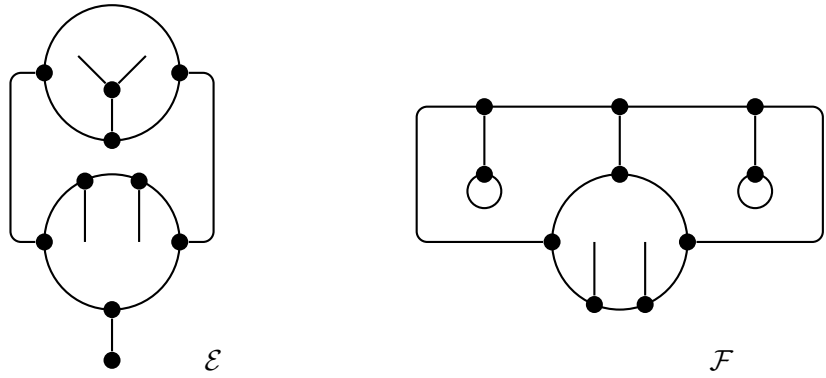

%%%%%%%%%%%%%%%%%%%%

\begin{figure}[h!]
\begin{center}
\begin{tikzpicture}[scale=0.45, inner sep=0.8mm]

\draw [thick] (2,5) arc (0:360:2);
\node (a1) at (2,5) [shape=circle, fill=black] {};
\node (a1') at (-2,5) [shape=circle, fill=black] {};
\node (b1) at (0,3) [shape=circle, fill=black] {};
\node (c1) at (0,4.5) [shape=circle, fill=black] {};

\draw [thick] (b1) to (c1);
\draw [thick] (1,5.5) to (c1) to (-1,5.5);

\draw [thick] (2,0) arc (0:360:2);
\node (a2) at (2,0) [shape=circle, fill=black] {};
\node (a2') at (-2,0) [shape=circle, fill=black] {};
\node (b2) at (0,-2) [shape=circle, fill=black] {};
\node (c2) at (0,-0.5) [shape=circle, fill=black] {};

\draw [thick] (b2) to (c2);
\draw [thick] (1,0.5) to (c2) to (-1,0.5);

\draw [thick] (2,-5) arc (0:360:2);
\node (a3) at (2,-5) [shape=circle, fill=black] {};
\node (a3') at (-2,-5) [shape=circle, fill=black] {};
\node (b3) at (0,-7) [shape=circle, fill=black] {};
\node (c3) at (0,-5.5) [shape=circle, fill=black] {};

\draw [thick] (b3) to (c3);
\draw [thick] (1,-4.5) to (c3) to (-1,-4.5);

\node (d) at (4,0) [shape=circle, fill=black] {};
\node (d') at (-4,0) [shape=circle, fill=black] {};
\draw [rounded corners, thick] (a1) to (4,5) to (d) to (4,-5) to (a3);
\draw [rounded corners, thick] (a1') to (-4,5) to (d') to (-4,-5) to (a3');
\draw [thick] (a2) to (d);
\draw [thick] (a2') to (d');

\node at (5, -7) {$\G$};

%%%%%%

\draw [thick] (17,5) arc (0:360:2);
\node (a1) at (17,5) [shape=circle, fill=black] {};
\node (a1') at (13,5) [shape=circle, fill=black] {};
\node (b1) at (15,3) [shape=circle, fill=black] {};
\node (c1) at (15,4.5) [shape=circle, fill=black] {};

\draw [thick] (b1) to (c1);
\draw [thick] (16,5.5) to (c1) to (14,5.5);

\draw [thick] (17,0) arc (0:360:2);
\node (a2) at (17,0) [shape=circle, fill=black] {};
\node (a2') at (13,0) [shape=circle, fill=black] {};
\node (c2) at (16,-1.75) [shape=circle, fill=black] {};
\node (c2') at (14,-1.75) [shape=circle, fill=black] {};
\node (d2) at (15,-2) [shape=circle, fill=black] {};

\draw [thick] (c2) to (16,0);
\draw [thick] (c2') to (14,0);

\node (e) at (19,0) [shape=circle, fill=black] {};
\node (e') at (11,0) [shape=circle, fill=black] {};
\node (f) at (17.5,-3.5) [shape=circle, fill=black] {};
\node (f') at (12.5,-3.5) [shape=circle, fill=black] {};
\node (g) at (15,-3.5) [shape=circle, fill=black] {};
\draw [thick] (g) to (d2);
\draw [rounded corners, thick] (a1) to (19,5) to (19,-3.5) to (11,-3.5) to (11,5) to (a1');
\draw [thick] (a2) to (e);
\draw [thick] (a2') to (e');
\draw [thick] (f) to (17.5,-2);
\draw [thick] (f') to (12.5,-2);

\node at (15, -7) {$\H$};

\end{tikzpicture}

\end{center}
\caption{The maps $\G$ and $\H$.}
\label{DessinsGH}
\end{figure}
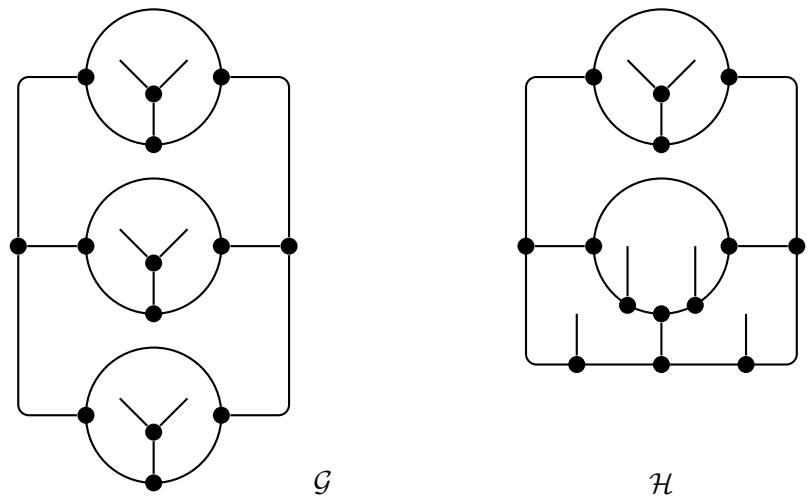

\newpage

%%%%%%%%%%%%%%%%%

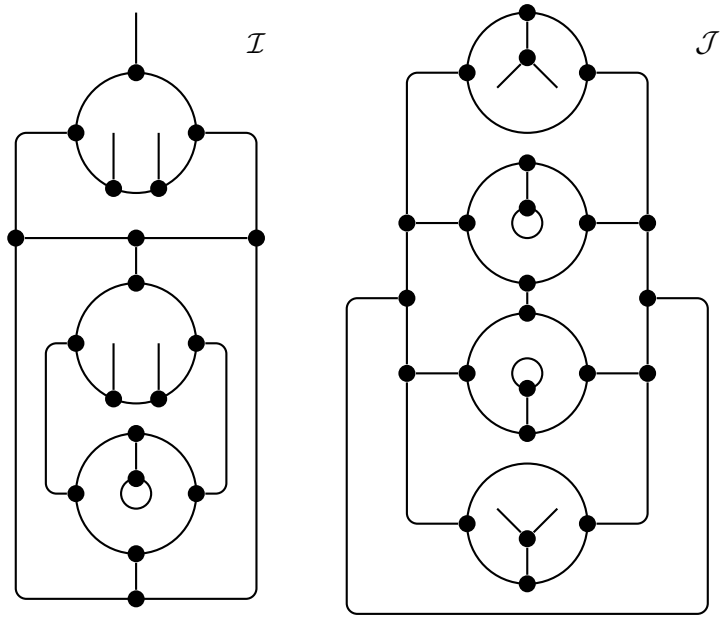
\begin{figure}[h!]
\begin{center}
\begin{tikzpicture}[scale=0.4, inner sep=0.8mm]

\draw [thick] (-11,8) arc (0:360:2);
\node (A1) at (-11,8) [shape=circle, fill=black] {};
\node (A1') at (-15,8) [shape=circle, fill=black] {};
\node (B1) at (-13,10) [shape=circle, fill=black] {};
\node (C1) at (-12.25,6.15) [shape=circle, fill=black] {};
\node (C1') at (-13.75,6.15) [shape=circle, fill=black] {};

\draw [thick] (C1) to (-12.25,8);
\draw [thick] (C1') to (-13.75,8);
\draw [thick] (B1) to (-13,12);

%%%%%

\draw [thick] (-11,1) arc (0:360:2);
\node (A2) at (-11,1) [shape=circle, fill=black] {};
\node (A2') at (-15,1) [shape=circle, fill=black] {};
\node (B2) at (-13,3) [shape=circle, fill=black] {};
\node (C2) at (-12.25,-0.85) [shape=circle, fill=black] {};
\node (C2') at (-13.75,-0.85) [shape=circle, fill=black] {};

\draw [thick] (C2) to (-12.25,1);
\draw [thick] (C2') to (-13.75,1);

%%%%%%

\draw [thick] (-11,-4) arc (0:360:2);
\node (A3) at (-11,-4) [shape=circle, fill=black] {};
\node (A3') at (-15,-4) [shape=circle, fill=black] {};
\node (B3) at (-13,-2) [shape=circle, fill=black] {};
\node (C3) at (-13,-3.5) [shape=circle, fill=black] {};
\node (D3) at (-13,-6) [shape=circle, fill=black] {};

\draw [thick] (B3) to (C3);
\draw [thick] (-12.5,-4) arc (0:360:0.5);

%%%%%

\node (E) at (-13,4.5) [shape=circle, fill=black] {};
\node (F) at (-9,4.5) [shape=circle, fill=black] {};
\node (F') at (-17,4.5) [shape=circle, fill=black] {};
\node (G) at (-13,-7.5) [shape=circle, fill=black] {};
\draw [thick] (F) to (F');
\draw [rounded corners, thick] (A1) to (-9,8) to (-9,-7.5) to (-17,-7.5) to (-17,8) to (A1');
\draw [thick] (E) to (B2);
\draw [thick] (G) to (D3);
\draw [rounded corners, thick] (A2) to (-10,1) to (-10,-4) to (A3);
\draw [rounded corners, thick] (A2') to (-16,1) to (-16,-4) to (A3');

\node at (-9,11) {$\I$};

%%%%%%%%
%%%%%%%%

\draw [thick] (2,10) arc (0:360:2);
\node (a1) at (2,10) [shape=circle, fill=black] {};
\node (a1') at (-2,10) [shape=circle, fill=black] {};
\node (b1) at (0,12) [shape=circle, fill=black] {};
\node (c1) at (0,10.5) [shape=circle, fill=black] {};

\draw [thick] (b1) to (c1);
\draw [thick] (1,9.5) to (c1) to (-1,9.5);

%%%%%%%

\draw [thick] (2,5) arc (0:360:2);
\node (a2) at (2,5) [shape=circle, fill=black] {};
\node (a2') at (-2,5) [shape=circle, fill=black] {};
\node (b2) at (0,7) [shape=circle, fill=black] {};
\node (c2) at (0,5.5) [shape=circle, fill=black] {};
\node (d2) at (0,3) [shape=circle, fill=black] {};

\draw [thick] (b2) to (c2);
\draw [thick] (0.5,5) arc (0:360:0.5);

%%%%%%%

\draw [thick] (2,0) arc (0:360:2);
\node (a3) at (2,0) [shape=circle, fill=black] {};
\node (a3') at (-2,0) [shape=circle, fill=black] {};
\node (b3) at (0,-2) [shape=circle, fill=black] {};
\node (c3) at (0,-0.5) [shape=circle, fill=black] {};
\node (d3) at (0,2) [shape=circle, fill=black] {};

\draw [thick] (b3) to (c3);
\draw [thick] (0.5,0) arc (0:360:0.5);

%%%%%%%

\draw [thick] (2,-5) arc (0:360:2);
\node (a4) at (2,-5) [shape=circle, fill=black] {};
\node (a4') at (-2,-5) [shape=circle, fill=black] {};
\node (b4) at (0,-7) [shape=circle, fill=black] {};
\node (c4) at (0,-5.5) [shape=circle, fill=black] {};

\draw [thick] (b4) to (c4);
\draw [thick] (1,-4.5) to (c4) to (-1,-4.5);

%%%%%%%

\node (e2) at (4,5) [shape=circle, fill=black] {};
\node (e2') at (-4,5) [shape=circle, fill=black] {};
\node (e3) at (4,0) [shape=circle, fill=black] {};
\node (e3') at (-4,0) [shape=circle, fill=black] {};
\draw [thick] (a2) to (e2);
\draw [thick] (a2') to (e2');
\draw [thick] (a3) to (e3);
\draw [thick] (a3') to (e3');
\draw [rounded corners, thick] (a1) to (4,10) to (e2) to (e3) to (4,-5) to (a4);
\draw [rounded corners, thick] (a1') to (-4,10) to (e2') to (e3') to (-4,-5) to (a4');

%%%%%%%%

\node (f) at (4,2.5) [shape=circle, fill=black] {};
\node (f') at (-4,2.5) [shape=circle, fill=black] {};
\draw [rounded corners, thick] (f) to (6,2.5) to (6,-8) to (-6,-8) to (-6,2.5) to (f');
\draw [thick] (d2) to (d3);

\node at (6,11) {$\J$};

\end{tikzpicture}

\end{center}
\caption{The maps $\I$ and $\J$.}
\label{DessinsIJ}
\end{figure}

%\newpage

%%%%%%%%%%%%%%%%%

\begin{figure}[h!]
\begin{center}
\begin{tikzpicture}[scale=0.4, inner sep=0.8mm]

\node (x) at (0,11) {};

\draw [thick] (2,7) arc (0:360:2);
\node (a) at (2,7) [shape=circle, fill=black] {};
\node (a') at (-2,7) [shape=circle, fill=black] {};
\node (b) at (0,5) [shape=circle, fill=black] {};
\node (c) at (0,6.5) [shape=circle, fill=black] {};

\draw [thick] (b) to (c);
\draw [thick] (1,7.5) to (c) to (-1,7.5);

%%%%%

\draw [thick] (-2,0) arc (0:360:2);
\node (a1) at (-2,0) [shape=circle, fill=black] {};
\node (a1') at (-6,0) [shape=circle, fill=black] {};
\node (b1) at (-4,-2) [shape=circle, fill=black] {};
\node (c1) at (-4,-0.5) [shape=circle, fill=black] {};
\node (d1) at (-4,2) [shape=circle, fill=black] {};

\draw [thick] (b1) to (c1);
\draw [thick] (-3.5,0) arc (0:360:0.5);

%%%%%%%
\draw [thick] (6,0) arc (0:360:2);
\node (a2) at (6,0) [shape=circle, fill=black] {};
\node (a2') at (2,0) [shape=circle, fill=black] {};
\node (b2) at (4,-2) [shape=circle, fill=black] {};
\node (c2) at (4,-0.5) [shape=circle, fill=black] {};
\node (d2) at (4,2) [shape=circle, fill=black] {};

\draw [thick] (b2) to (c2);
\draw [thick] (4.5,0) arc (0:360:0.5);

%%%%%

\node (e) at (4,3.5) [shape=circle, fill=black] {};
\node (f) at (9,3.5) [shape=circle, fill=black] {};
\node (g) at (9,0) [shape=circle, fill=black] {};
\node (h) at (7.5,0) [shape=circle, fill=black] {};
\node (e') at (-4,3.5) [shape=circle, fill=black] {};
\node (f') at (-9,3.5) [shape=circle, fill=black] {};
\node (g') at (-9,0) [shape=circle, fill=black] {};
\node (h') at (-7.5,0) [shape=circle, fill=black] {};
\node (i) at (0,-3.5) [shape=circle, fill=black] {};
\node (j) at (0,0) [shape=circle, fill=black] {};

\draw [thick] (f) to (f');
\draw [thick] (g) to (a2);
\draw [thick] (g') to (a1');
\draw [thick] (a1) to (a2');
\draw [thick] (e) to (d2);
\draw [thick] (e') to (d1);
\draw [thick] (i) to (j);
\draw [thick] (h) to (7.5,1.5);
\draw [thick] (h') to (-7.5,1.5);
\draw [rounded corners, thick] (a') to (-9,7) to (-9,-3.5) to (9,-3.5) to (9,7) to (a);

\end{tikzpicture}

\end{center}
\caption{The map $\K$.}
\label{dessinK}
\end{figure}

%%%%%%%%%%%%%%%%%%%%%%%%%

\begin{figure}[h!]
\begin{center}
\begin{tikzpicture}[scale=0.4, inner sep=0.8mm]

\draw [thick] (13,0) arc (0:360:2);
\node (a1) at (13,0) [shape=circle, fill=black] {};
\node (a'1) at (9,0) [shape=circle, fill=black] {};
\node (b1) at (11,-2) [shape=circle, fill=black] {};
\node (c1) at (11,-0.5) [shape=circle, fill=black] {};
\node (d1) at (11,2) [shape=circle, fill=black] {};

\draw [thick] (b1) to (c1);
\draw [thick] (11.5,0) arc (0:360:0.5);

%%%%%

\draw [thick] (-9,0) arc (0:360:2);
\node (a2) at (-9,0) [shape=circle, fill=black] {};
\node (a'2) at (-13,0) [shape=circle, fill=black] {};
\node (b2) at (-11,-2) [shape=circle, fill=black] {};
\node (c2) at (-11,-0.5) [shape=circle, fill=black] {};
\node (d2) at (-11,2) [shape=circle, fill=black] {};

\draw [thick] (b2) to (c2);
\draw [thick] (-10.5,0) arc (0:360:0.5);

%%%%%%%

\draw [thick] (2,0) arc (0:360:2);
\node (a) at (2,0) [shape=circle, fill=black] {};
\node (a') at (-2,0) [shape=circle, fill=black] {};
\node (b) at (0,2) [shape=circle, fill=black] {};
\node (c) at (0.75,-1.85) [shape=circle, fill=black] {};
\node (c') at (-0.75,-1.85) [shape=circle, fill=black] {};

\draw [thick] (c) to (0.75,0);
\draw [thick] (c') to (-0.75,0);

%%%%%%%

\draw [thick] (6,-8) arc (0:360:2);
\node (a3) at (6,-8) [shape=circle, fill=black] {};
\node (a'3) at (2,-8) [shape=circle, fill=black] {};
\node (b3) at (4,-6) [shape=circle, fill=black] {};
\node (c3) at (4,-7.5) [shape=circle, fill=black] {};
\node (d3) at (4,-10) [shape=circle, fill=black] {};

\draw [thick] (b3) to (c3);
\draw [thick] (4.5,-8) arc (0:360:0.5);

%%%%%%%%

\draw [thick] (-2,-8) arc (0:360:2);
\node (a4) at (-2,-8) [shape=circle, fill=black] {};
\node (a'4) at (-6,-8) [shape=circle, fill=black] {};
\node (b4) at (-4,-6) [shape=circle, fill=black] {};
\node (c4) at (-4,-7.5) [shape=circle, fill=black] {};
\node (d4) at (-4,-10) [shape=circle, fill=black] {};

\draw [thick] (b4) to (c4);
\draw [thick] (-3.5,-8) arc (0:360:0.5);

%%%%%%%%

\draw [thick] (a) to (a'1);
\draw [thick] (a') to (a2);

\node (e) at (4,0) [shape=circle, fill=black] {};
\node (e') at (-4,0) [shape=circle, fill=black] {};
\draw [rounded corners, thick] (e) to (4,-4) to (-4,-4) to (e');

\node (f) at (0,-4) [shape=circle, fill=black] {};
\node (g) at (0,-8) [shape=circle, fill=black] {};
\draw [thick] (f) to (g);
\draw [thick] (a'3) to (g) to (a4);

\node (h) at (7.5,0) [shape=circle, fill=black] {};
\node (h') at (-7.5,0) [shape=circle, fill=black] {};
\draw [rounded corners, thick] (h) to (7.5,-8) to (a3);
\draw [rounded corners, thick] (h') to (-7.5,-8) to (a'4);

\node (i) at (4,-12) [shape=circle, fill=black] {};
\node (i') at (-4,-12) [shape=circle, fill=black] {};
\draw [rounded corners, thick] (a1) to (14.5,0) to (14.5,-12) to (i) to (i') to (-14.5,-12) to (-14.5,0) to (a'2);
\draw [thick] (i) to (d3);
\draw [thick] (i') to (d4);

\node (j) at (0,4) [shape=circle, fill=black] {};
\draw [thick] (j) to (b);
\draw [rounded corners, thick] (d1) to (11,4) to (j) to (-11,4) to (d2);

\end{tikzpicture}

\end{center}
\caption{The map $\L$.}
\label{dessinL}
\end{figure}

\newpage

%%%%%%%%%%%%%%%%%%%%%%%

\begin{figure}[h!]
\begin{center}
\begin{tikzpicture}[scale=0.4, inner sep=0.8mm]

\draw [thick] (10,0) arc (0:360:2);
\node (a1) at (10,0) [shape=circle, fill=black] {};
\node (a'1) at (6,0) [shape=circle, fill=black] {};
\node (b1) at (8,-2) [shape=circle, fill=black] {};
\node (c1) at (8,-0.5) [shape=circle, fill=black] {};
\node (d1) at (8,2) [shape=circle, fill=black] {};

\draw [thick] (b1) to (c1);
\draw [thick] (8.5,0) arc (0:360:0.5);

%%%%%%%

\draw [thick] (-6,0) arc (0:360:2);
\node (a2) at (-6,0) [shape=circle, fill=black] {};
\node (a'2) at (-10,0) [shape=circle, fill=black] {};
\node (b2) at (-8,-2) [shape=circle, fill=black] {};
\node (c2) at (-8,-0.5) [shape=circle, fill=black] {};
\node (d2) at (-8,2) [shape=circle, fill=black] {};

\draw [thick] (b2) to (c2);
\draw [thick] (-7.5,0) arc (0:360:0.5);

%%%%%%%

\draw [thick] (2,-3) arc (0:360:2);
\node (a3) at (2,-3) [shape=circle, fill=black] {};
\node (a'3) at (-2,-3) [shape=circle, fill=black] {};
\node (b3) at (0,-5) [shape=circle, fill=black] {};
\node (c3) at (0,-3.5) [shape=circle, fill=black] {};

\draw [thick] (b3) to (c3);
\draw [thick] (1,-2.5) to (c3) to (-1,-2.5);

%%%%%%%

\draw [thick] (2,4) arc (0:360:2);
\node (a4) at (2,4) [shape=circle, fill=black] {};
\node (a'4) at (-2,4) [shape=circle, fill=black] {};
\node (b4) at (0,2) [shape=circle, fill=black] {};
\node (c4) at (0.75,5.85) [shape=circle, fill=black] {};
\node (c'4) at (-0.75,5.85) [shape=circle, fill=black] {};

\draw [thick] (c4) to (0.75,4);
\draw [thick] (c'4) to (-0.75,4);

%%%%%%%

\draw [thick] (2,11) arc (0:360:2);
\node (a5) at (2,11) [shape=circle, fill=black] {};
\node (a'5) at (-2,11) [shape=circle, fill=black] {};
\node (b5) at (0,13) [shape=circle, fill=black] {};
\node (c5) at (0,11.5) [shape=circle, fill=black] {};
\node (d5) at (0,9) [shape=circle, fill=black] {};

\draw [thick] (b5) to (c5);
\draw [thick] (0.5,11) arc (0:360:0.5);

%%%%%%%%%

\draw [thick] (a'1) to (a2);
\node (e) at (0,0) [shape=circle, fill=black] {};
\node (f) at (4,0) [shape=circle, fill=black] {};
\node (f') at (-4,0) [shape=circle, fill=black] {};
\node (g) at (4,-3) [shape=circle, fill=black] {};
\node (g') at (-4,-3) [shape=circle, fill=black] {};
\node (h) at (4,-6) [shape=circle, fill=black] {};
\node (h') at (-4,-6) [shape=circle, fill=black] {};
\draw [thick] (e) to (b4);
\draw [thick] (f) to (h);
\draw [thick] (f') to (h');
\draw [thick] (a3) to (g);
\draw [thick] (a'3) to (g');

%%%%%%

\node (i) at (12,0) [shape=circle, fill=black] {};
\node (i') at (-12,0) [shape=circle, fill=black] {};
\draw [thick] (a1) to (i);
\draw [thick] (a'2) to (i');
\draw [rounded corners, thick] (a5) to (12,11) to (12,-6) to (-12,-6) to (-12,11) to (a'5);

%%%%%%%

\node (j) at (0,7.5) [shape=circle, fill=black] {};
\node (k) at (8,4) [shape=circle, fill=black] {};
\node (k') at (-8,4) [shape=circle, fill=black] {};
\draw [rounded corners, thick] (d1) to (8,7.5) to (-8,7.5) to (d2);
\draw [thick] (k) to (a4);
\draw [thick] (k') to (a'4);
\draw [thick] (j) to (d5);

\end{tikzpicture}

\end{center}
\caption{The map $\M$.}
\label{dessinM}
\end{figure}

%%%%%%%%%%%%%%%%%%

\begin{figure}[h!]
\begin{center}
\begin{tikzpicture}[scale=0.4, inner sep=0.8mm]

\draw [thick] (2,0) arc (0:360:2);
\node (a1) at (2,0) [shape=circle, fill=black] {};
\node (a'1) at (-2,0) [shape=circle, fill=black] {};
\node (b1) at (0,-2) [shape=circle, fill=black] {};
\node (c1) at (1,-1.7) [shape=circle, fill=black] {};
\node (c'1) at (-1,-1.7) [shape=circle, fill=black] {};

\draw [thick] (c1) to (1,0);
\draw [thick] (c'1) to (-1,0);

%%%%%%%

\draw [thick] (2,5) arc (0:360:2);
\node (a2) at (2,5) [shape=circle, fill=black] {};
\node (a'2) at (-2,5) [shape=circle, fill=black] {};
\node (b2) at (0,3) [shape=circle, fill=black] {};
\node (c2) at (0,4.5) [shape=circle, fill=black] {};
\node (d2) at (0,7) [shape=circle, fill=black] {};

\draw [thick] (b2) to (c2);
\draw [thick] (0.5,5) arc (0:360:0.5);

%%%%%%%%

\draw [thick] (12,0) arc (0:360:2);
\node (a3) at (12,0) [shape=circle, fill=black] {};
\node (a'3) at (8,0) [shape=circle, fill=black] {};
\node (b3) at (10,-2) [shape=circle, fill=black] {};
\node (c3) at (10,-0.5) [shape=circle, fill=black] {};
\node (d3) at (10,2) [shape=circle, fill=black] {};

\draw [thick] (b3) to (c3);
\draw [thick] (10.5,0) arc (0:360:0.5);

%%%%%%%%

\draw [thick] (-8,0) arc (0:360:2);
\node (a4) at (-8,0) [shape=circle, fill=black] {};
\node (a'4) at (-12,0) [shape=circle, fill=black] {};
\node (b4) at (-10,-2) [shape=circle, fill=black] {};
\node (c4) at (-10,-0.5) [shape=circle, fill=black] {};
\node (d4) at (-10,2) [shape=circle, fill=black] {};

\draw [thick] (b4) to (c4);
\draw [thick] (-9.5,0) arc (0:360:0.5);

%%%%%%%%

\draw [thick] (2,-8) arc (0:360:2);
\node (a5) at (2,-8) [shape=circle, fill=black] {};
\node (a'5) at (-2,-8) [shape=circle, fill=black] {};
\node (b5) at (0,-10) [shape=circle, fill=black] {};
\node (c5) at (0,-8.5) [shape=circle, fill=black] {};

\draw [thick] (b5) to (c5);
\draw [thick] (1,-7.5) to (c5) to (-1,-7.5);

%%%%%%%

\draw [thick] (a1) to (a'3);
\draw [thick] (a'1) to (a4);
\node (e) at (4,0) [shape=circle, fill=black] {};
\node (e') at (-4,0) [shape=circle, fill=black] {};
\node (f) at (6,0) [shape=circle, fill=black] {};
\node (f') at (-6,0) [shape=circle, fill=black] {};
\node (g) at (6,-4) [shape=circle, fill=black] {};
\node (g') at (-6,-4) [shape=circle, fill=black] {};
\draw [thick] (f) to (g);
\draw [thick] (f') to (g');
\node (h) at (0,-4) [shape=circle, fill=black] {};
\node (i) at (13.5,-4) [shape=circle, fill=black] {};
\node (i') at (-13.5,-4) [shape=circle, fill=black] {};
\draw [thick] (i) to (i');
\draw [thick] (h) to (b1);

%%%%%%%

\draw [rounded corners, thick] (a2) to (4,5) to (e);
\draw [rounded corners, thick] (a'2) to (-4,5) to (e');
\node (j) at (0,9) [shape=circle, fill=black] {};
\draw [thick] (j) to (d2);
\draw [rounded corners, thick] (d3) to (10,9) to (-10,9) to (d4);
\node (k) at (13.5,0) [shape=circle, fill=black] {};
\node (k') at (-13.5,0) [shape=circle, fill=black] {};
\draw [thick] (a3) to (k);
\draw [thick] (a'4) to (k');
\draw [rounded corners, thick] (a5) to (13.5,-8) to (13.5,10) to (-13.5,10) to (-13.5,-8) to (a'5);

\end{tikzpicture}

\end{center}
\caption{The map $\N$.}
\label{dessinN}
\end{figure}

\bigskip\bigskip

{
\noindent {\sc G. A. Jones:}\\ School of Mathematics,\\ University of Southampton,\\ Southampton SO17 1BJ,\\ U.K.\\ \medskip email: {\tt G.A.Jones@maths.soton.ac.uk}

\bigskip

\noindent {\sc E. Pierro:}\\ Department of Economics, Mathematics and Statistics,\\ Birkbeck, University of London,\\ Malet Street,\\ London WC1E 7HX,\\ U.K.\\ \medskip email: {\tt e.pierro@mail.bbk.ac.uk}
}

\end{document}